\RequirePackage{fix-cm}
\documentclass[smallextended]{svjour3}       
\smartqed  
\usepackage{graphicx}
\usepackage{algorithm}
\usepackage{algorithmic}

\usepackage{float}
\usepackage{amsfonts}
\usepackage{verbatim}
\usepackage{amssymb}
\begin{document}

\title{Exact recovery low-rank matrix via transformed affine matrix rank minimization}

\author{Angang Cui$^{1}$\and
        Jigen Peng$^{2,\ast}$\and
        Haiyang Li$^{3}$
}


\institute{$\ast$ Corresponding author\\
            Jigen Peng \at
            \email{jgpengxjtu@126.com}\\
             1 School of Mathematics and Statistics, Xi'an Jiaotong University, Xi'an, 710049, China;\\
             2 School of Mathematics and Information Science, Guangzhou University, Guangzhou, 510006, China;\\
             3 School of Science, Xi'an Polytechnic University, Xi'an, 710048, China.
}

\date{Received: date / Accepted: date}

\maketitle

\begin{abstract}
The goal of affine matrix rank minimization problem is to reconstruct a low-rank or approximately low-rank matrix under linear constraints. In general,
this problem is combinatorial and NP-hard. In this paper, a nonconvex fraction function is studied to approximate the rank of a matrix and translate
this NP-hard problem into a transformed affine matrix rank minimization problem. The equivalence between these two problems is established, and we
proved that the uniqueness of the global minimizer of transformed affine matrix rank minimization problem also solves affine matrix rank minimization
problem if some conditions are satisfied. Moreover, we also proved that the optimal solution to the transformed affine matrix rank minimization problem can
be approximately obtained by solving its regularization problem for some proper smaller $\lambda>0$. Lastly, the DC algorithm is utilized to solve
the regularization transformed affine matrix rank minimization problem and the numerical experiments on image inpainting problems show that our method
performs effectively in recovering low-rank images compared with some state-of-art algorithms.

\end{abstract}

\section{Introduction}\label{section1}
The goal of affine matrix rank minimization (AMRM) problem is to reconstruct a low-rank or approximately low-rank matrix that satisfies a given system of linear equality constraints.
In mathematics, it can be described as the following minimization problem
\begin{equation}\label{equ1}
(\mathrm{AMRM})\ \ \ \ \ \ \min_{X\in \mathbb{R}^{m\times n}} \ \mbox{rank}(X)\ \ s.t. \ \  \mathcal{A}(X)=b
\end{equation}
where $\mathcal{A}: \mathbb{R}^{m\times n}\mapsto \mathbb{R}^{d}$ is the linear map and the vector $b\in \mathbb{R}^{d}$. Without loss of generality, we assume $m\leq n$. Many
applications arising in various areas can be captured by solving the problem (AMRM), for instance, the network localization \cite{Ji1}, the minimum order system and low-dimensional
Euclidean embedding in control theory \cite{Faz2,Faze3}, the collaborative filtering in recommender systems \cite{Cand4,Jan5}, and so on.
One important special case of the problem (AMRM) is the matrix completion (MC) problem \cite{Cand4}
\begin{equation}\label{equ2}
(\mathrm{MC})\ \ \ \ \ \ \min_{X\in \mathbb{R}^{m\times n}} \ \mbox{rank}(X)\ \ s.t. \ \  X_{i,j}=M_{i,j},\ \ (i,j)\in \Omega.
\end{equation}
This completion problem has been applied in the famous Netflix problem \cite{Netfix6}, image inpainting problem \cite{Yeg7} and machine learning \cite{Xi8,Xi9}.
In general, however, the problem (AMRM) is a challenging non-convex optimization problem and is known as NP-hard \cite{Recht10} due to the combinational nature of the rank function.

Among the numerous substitution models, the nuclear-norm affine matrix rank minimization (NAMRM) problem has been considered as the most popular alternative
\cite{Faze3,Cand4,Cand11,Faze12,Cand13}:
\begin{equation}\label{equ3}
(\mathrm{NAMRM})\ \ \ \ \ \ \min_{X\in \mathbb{R}^{m\times n}} \ \|X\|_{\ast}\ \ s.t. \ \  \mathcal{A}(X)=b.
\end{equation}
where $\|X\|_{\ast}=\sum_{i=1}^{m}\sigma_{i}(X)$ is the nuclear-norm of the matrix $X\in \mathbb{R}^{m\times n}$. Recht et al. in \cite{Recht10} have show that if a certain restricted
isometry property (RIP) holds for the linear transformation defining the constraints, the minimum rank solution of problem (AMRM) can be recovered by solving the problem (NAMRM).
In addition, some popular methods, including singular value thresholding algorithm \cite{Cai14}, proximal gradient algorithm \cite{Ma15} and accelerated  proximal gradient
algorithm \cite{Toh16}, are proposed to solve its regularization (or Lagrangian) version:
\begin{equation}\label{equ4}
(\mathrm{RNAMRM})\ \ \ \ \ \ \min_{X\in \mathbb{R}^{m\times n}} \Big\{\|\mathcal{A}(X)-b\|_{2}^{2}+\lambda\|X\|_{\ast}\Big\}
\end{equation}
where $\lambda>0$ is the regularization parameter can be selected to guarantee that solutions of the problem (NAMRM) and (RNAMRM) are same \cite{Yu17}. However, these algorithms tend to have
biased estimation by shrinking all the singular values toward zero simultaneously, and sometimes results in over-penalization in the regularization problem (RNAMRM) as the $\ell_{1}$-norm in
compressive sensing. Moreover, with the recent development of non-convex relaxation approach in sparse signal recovery problems, many researchers have shown that using a non-convex surrogate
function to approximate the $\ell_{0}$-norm is a better choice than using the $\ell_{1}$-norm. This brings our attention back to the non-convex surrogate functions of the rank function.

In this paper, a continuous promoting low-rank non-convex function
\begin{equation}\label{equ5}
P_{a}(X)=\sum_{i=1}^{m}\rho_{a}(\sigma_{i}(X))=\sum_{i=1}^{m}\frac{a\sigma_{i}(X)}{a\sigma_{i}(X)+1}
\end{equation}
in terms of the singular values of matrix $X$ is considered to substitute the rank function $\mathrm{rank}(X)$ in the problem (AMRM), where the non-convex function
\begin{equation}\label{equ6}
\rho_{a}(t)=\frac{a|t|}{a|t|+1}\ \ \ (a>0)
\end{equation}
is the fraction function. It is to see clearly that, with the change of parameter $a>0$, the non-convex function $P_{a}(X)$ approximates the rank of matrix $X$:
\begin{equation}\label{equ7}
\lim_{a\rightarrow+\infty}P_{a}(X)=\lim_{a\rightarrow+\infty}\sum_{i=1}^{m}\frac{a\sigma_{i}(X)}{a\sigma_{i}(X)+1}\approx\left\{
    \begin{array}{ll}
      0, & {\mathrm{if} \ \sigma_{i}(X)=0;} \\
      \mathrm{rank}(X), & {\mathrm{if} \ \sigma_{i}(X)> 0.}
    \end{array}
  \right.
\end{equation}
By this transformation, the NP-hard problem (AMRM) could be relaxed into the following matrix rank minimization problem with a continuous non-convex penalty, namely,
transformed affine matrix rank minimization (TrAMRM) problem:
\begin{equation}\label{equ8}
(\mathrm{TrAMRM})\ \ \ \ \ \ \ \ \min_{X\in \mathbb{R}^{m\times n}} \ P_{a}(X)\ \ s.t. \ \  \mathcal{A}(X)=b
\end{equation}
where the non-convex surrogate function $P_{a}(X)$ in terms of the singular values of matrix $X$ is defined in (\ref{equ5}). Unfortunately, although we relax the NP-hard
problem (AMRM) into a continuous problem (TrAMRM), this relax problem is still computationally harder to solve due to the non-convex nature of the function $P_{a}(X)$, in
fact it is also NP-hard. Frequently, we consider its regularization version:
\begin{equation}\label{equ9}
(\mathrm{RTrAMRM})\ \ \ \ \ \min_{X\in \mathbb{R}^{m\times n}} \Big\{\|\mathcal{A}(X)-b\|_{2}^{2}+\lambda P_{a}(X)\Big\}
\end{equation}
where $\lambda>0$ is the regularization parameter. Unlike the convex optimal theory, there are no parameters $\lambda>0$ such that the solution to the regularization
problem (RTrAMRM) also solves the constrained problem (TrAMRM). However, as the unconstrained form, the problem (RNuAMRM) may possess much more algorithmic advantages.
Moreover, we also proved that the optimal solution to the problem (TrAMRM) can be approximately obtained by solving the problem (RTrAMRM) for some proper smaller $\lambda>0$.

The rest of this paper is organized as follows. Some notions and preliminary results that are used in this paper are given in Section \ref{section2}.
In Section \ref{section3}, the equivalence of the problem (TrAMRM) and (AMRM) is established. Moreover, we proved that the optimal solution to the
problem (TrAMRM) can be approximately obtained by solving the problem (RTrAMRM) for some proper smaller $\lambda>0$. In Section \ref{section4}, the
DC algorithm is utilized to solve the problem (RTrAMRM) and the numerical results of the numerical experiments on image inpainting problems are
demonstrated in Section \ref{section5}. Finally, we give some concluding remarks in Section \ref{section6}.

\section{Preliminaries} \label{section2}
In this section, we give some notions and preliminary results that are used in this paper.

\subsection{Some notions} \label{subsection2.1}
The linear map $\mathcal{A}: \mathbb{R}^{m\times n}\mapsto \mathbb{R}^{d}$ determined by $d$ matrices $A_{1}, A_{2}, \cdots, A_{d}\in \mathbb{R}^{m\times n}$ can
be expressed as $\mathcal{A}(X)=(\langle A_{1},X\rangle, \langle A_{2},X\rangle,\cdots, \langle A_{d},X\rangle)^{\top}\in \mathbb{R}^{d}$. Let $A=(vec(A_{1}),
vec(A_{2}),\cdots, vec(A_{d}))^{\top}\in \mathbb{R}^{d\times mn}$ and $x=vec(X)\in \mathbb{R}^{mn}$. Then we can get that $\mathcal{A}(X)=Ax$. The standard inner
product of matrices $X\in \mathbb{R}^{m\times n}$ and $Y\in \mathbb{R}^{m\times n}$ is denoted by $\langle X,Y\rangle$, and $\langle X,Y\rangle=\mathrm{tr}(Y^{\top}X)$.
The $\mathcal{A}^{\ast}$ denotes the adjoint of $\mathcal{A}$, and for any $y\in \mathbb{R}^{d}$, $\mathcal{A}^{\ast}(y)=\sum_{i=1}^{d}y_{i}A_{i}$. The singular value
decomposition (SVD) of matrix $X$ is $X=U\Sigma V^{\top}$, where $U$ is an $m\times m$  unitary matrix, $V$ is an $n\times n$ unitary matrix and
$\Sigma=\mathrm{Diag}(\sigma(X))\in \mathbb{R}^{m\times n}$ is a diagonal matrix. The vector $\sigma(X): \sigma_{1}(X)\geq\sigma_{2}(X)\geq\cdots\geq\sigma_{m}(X)$ arranged
in descending order represents the singular values vector of matrix $X$, and $\sigma_{i}(X)$ denotes the $i$-th largest singular value of matrix $X$ for $i=1,2,\cdots, m$.

\subsection{Some useful results} \label{subsection2.2}
\begin{lemma}\label{lem1}{\rm(see \cite{Recht10})}
Let $M$ and $N$ be matrices of the same dimensions. Then there exist matrices $N_{1}$ and $N_{2}$ such
that\\
(1) $N=N_{1}+N_{2}$;\\
(2) $\mathrm{rank}(N_{1})\leq 2\mathrm{rank}(M)$;\\
(3) $MN_{2}^{\top}=0$ and $M^{\top}N_{2}=0$;\\
(4) $\langle N1_{1}, N_{2}\rangle=0$.
\end{lemma}

By Lemma \ref{lem1}, we can derive the following important corollary.
\begin{corollary}\label{cor1}
Let $X^{\ast}$ and $X_{0}$ be the optimal solutions to the problem (TrAMRM) and (AMRM), respectively. If we set $R=X^{\ast}-X_{0}$, then there exist matrices $R_{0}$
and $R_{c}$ such that \\
(1) $R=R_{0}+R_{c}$;\\
(2) $\mathrm{rank}(R_{0})\leq 2\mathrm{rank}(X_{0})$;\\
(3) $X_{0}R_{c}^{\top}=0$, $X_{0}^{\top}R_{c}=0$ and $\langle R_{0}, R_{c}\rangle=0$.
\end{corollary}
\begin{lemma}\label{lem2}
Let $M$ and $N$ be matrices of the same dimensions. If $MN^{\top}=0$ and $M^{\top}N=0$, then $P_{a}(M+N)=P_{a}(M)+P_{a}(N)$.
\end{lemma}
\textbf{proof.} Consider the SVDs of the matrices $M$ and $N$:
\begin{equation}\label{equ10}
M=U_{M}\left[
         \begin{array}{cc}
           \Sigma_{M} & 0 \\
           0 & 0 \\
         \end{array}
       \right]V_{M}^{\top},\ \ \
N=U_{N}\left[
         \begin{array}{cc}
           \Sigma_{N} & 0 \\
           0 & 0 \\
         \end{array}
       \right]V_{N}^{\top}.
\end{equation}
Since $U_{M}$ and $U_{N}$ are left-invertible, the condition $MN^{\top}=0$ implies that $V_{M}^{\top}V_{N}=0$. Similarly,
$M^{\top}N=0$ implies that $U_{M }^{\top}U_{N}=0$. Thus, the following is a valid SVD for $M+N$,
\begin{equation}\label{equ11}
M+N=\left[
      \begin{array}{cc}
        U_{M} & U_{N} \\
      \end{array}
    \right]
\left[
  \begin{array}{cccc}
    \Sigma_{M} & 0 & 0 & 0 \\
    0 & 0 & 0 & 0 \\
    0 & 0 & \Sigma_{N} & 0 \\
    0 & 0 & 0 & 0 \\
  \end{array}
\right]
\left[
  \begin{array}{cc}
    V_{M} & V_{N} \\
  \end{array}
\right]^{\top}.
\end{equation}
This shows that the singular values of $M+N$ are equal to the union (with repetition) of the singular values of $M$ and $N$.
Hence, $P_{a}(M+N)=P_{a}(M)+P_{a}(N)$. This completes the proof. $\hfill{} \Box$

Combing Corollary \ref{cor1} and Lemma \ref{lem2}, we can get the following corollary.
\begin{corollary}\label{cor2}
Let $X^{\ast}$ and $X_{0}$ be the optimal solutions to the problem (TrAMRM) and (AMRM), respectively. If we set $R=X^{\ast}-X_{0}$, then, there exist matrices $R_{0}$
and $R_{c}$ such that $R=R_{0}+R_{c}$ and
\begin{equation}\label{equ12}
P_{a}(R_{c})\leq P_{a}(R_{0}).
\end{equation}
\end{corollary}
\textbf{proof.}
By optimality of $X^{\ast}$, we have $P_{a}(X_{0})\geq P_{a}(X^{\ast})$. Let $R=X^{\ast}-X_{0}$. Applying Corollary \ref{cor1} to the matrices $X_{0}$ and $R$, there
exist matrices $R_{0}$ and $R_{c}$ such that $R=R_{0}+R_{c}$, $\mathrm{rank}(R_{0})\leq 2\mathrm{rank}(X_{0})$, $X_{0}R_{c}^{\top}=0$, $X_{0}^{\top}R_{c}=0$. Then
\begin{equation}\label{equ13}
\begin{array}{llll}
P_{a}(X_{0})&\geq& P_{a}(X^{\ast})\\
&=& P_{a}(X_{0}+R)\\
&\geq& P_{a}(X_{0}+R_{c})-P_{a}(R_{0})\\
&=& P_{a}(X_{0})+P_{a}(R_{c})-P_{a}(R_{0})
\end{array}
\end{equation}
where the third assertion follows the triangle inequality and the last one follows Lemma \ref{lem2}. Rearranging (\ref{equ13}), we can conclude that
$$P_{a}(R_{c})\leq P_{a}(R_{0}).$$
This completes the proof.  $\hfill{} \Box$
\begin{definition}\label{def1}
Let $R=U \mathrm{Diag}(\sigma_{i}(R))V^{\top}$ be the SVD of the matrix $R$ defined in Corollary \ref{cor1}, we define the matrices $R_{0}$ and $R_{c}$ as:
\begin{equation}\label{equ14}
R_{0}=[U_{2T}\ 0]_{m\times m}\left[
                     \begin{array}{ccccccccc}
                       \Sigma_{1}  & & 0\\
                        & 0  &0
                     \end{array}
                   \right]_{m\times n}[V_{2T}\ 0]_{n\times n}^{\top}
\end{equation}
and
\begin{equation}\label{equ15}
R_{c}=[0\ U_{m-2T}]_{m\times m}\left[
                     \begin{array}{ccccccccc}
                       0 & & 0\\
                        & \Sigma_{2}&0
                     \end{array}
                   \right]_{m\times n}[0\ V_{n-2T}]_{n\times n}^{\top},
\end{equation}
where
$$\Sigma_{1}=\left[
               \begin{array}{ccc}
                 \sigma_{1}(R) &  &  \\
                  & \ddots &  \\
                  &  & \sigma_{2T}(R) \\
               \end{array}
             \right],\ \ \
\Sigma_{2}=\left[
               \begin{array}{ccc}
                 \sigma_{2T+1}(R) &  &  \\
                  & \ddots &  \\
                  &  & \sigma_{m}(R) \\
               \end{array}
             \right]
$$
and
$$U=[U_{2T},U_{m-2T}],\ \ \ V=[V_{2T},V_{m-2T}].$$
\end{definition}
\begin{definition}\label{def2}
For each positive integer $i\geq1$, we define the index set $I_{i}=\{K(i-1)+2T+1, \cdots, Ki+2T\}$ and partition matrix $R_{c}$ into a sum of matrices $R_{1},R_{2},\cdots$, i.e.,
$$R_{c}=\sum_{i}R_{i},$$
where
$$R_{i}=[0\ U_{I_{i}}\ 0]_{m\times m}\left[
                     \begin{array}{ccccccccc}
                       0 & & & & & & & &0\\
                        & \ddots & & & & & & &\vdots\\
                        &  & 0 & & & & & &0\\
                        &  &  & \sigma_{I_{i}} & & & & &0\\
                        &  &  & & 0 &  &  &  &0\\
                        & & & & &  &\ddots& &\vdots\\
                       & &  &  &  &  & &0 &0\\
                     \end{array}
                   \right]_{m\times n}[0\ V_{I_{i}}\ 0]_{n\times n}^{\top}.
$$
It is clear that $R_{i}^{\top}R_{h}=0$, $R_{i}R_{h}^{\top}=0$ for any $i\neq h$, and $\mathrm{rank}(R_{i})\leq K$.
\end{definition}

By the above lemmas and definitions, we shall derive some important results in this paper.
\begin{theorem}\label{the1}
The matrices $R_{0}\in \mathbb{R}^{m\times n}$ and $R_{1}\in \mathbb{R}^{m\times n}$ defined in Definition \ref{def1} and Definition \ref{def2} satisfy
\begin{equation}\label{equ16}
\|R_{0}+R_{1}\|_{F}\geq\frac{P_{a}(R_{0})}{a\sqrt{2T}}.
\end{equation}
\end{theorem}
\textbf{proof.}
Since
$$\rho_{a}(t)=\frac{a|t|}{a|t|+1}\leq a|t|,$$
we have
\begin{eqnarray*}
P_{a}(R_{0})&=&\displaystyle\sum_{i}\frac{a\sigma_{i}(R_{0})}{a\sigma_{i}(R_{0})+1}\\
&\leq& a\|\sigma(R_{0})\|_{1}\\
&\leq&a\sqrt{2T}\|R_{0}\|_{F}\\
&\leq& a\sqrt{2T}\|R_{0}+R_{1}\|_{F}.
\end{eqnarray*}
This completes the proof. $\hfill{} \Box$

\begin{theorem}\label{the2}
For any $\gamma>\displaystyle\frac{a(ar_{c}-a+1)\sigma_{1}(R_{c})}{a-1}$, the matrices $R_{i}$s defined in Definition \ref{def2} satisfy
\begin{equation}\label{equ17}
\sum_{i\geq2}\|\gamma^{-1}R_{i}\|_{F}\leq\sum_{i\geq2}\frac{P_{a}(\gamma^{-1}R_{i-1})}{\sqrt{K}}\leq\frac{P_{a}(\gamma^{-1}R_{0})}{\sqrt{K}}
\end{equation}
where $r_{c}=\mathrm{rank}(R_{c})$.
\end{theorem}

We will need the following technical lemma that shows for any matrix $X\in \mathbb{R}^{m\times n}$ there corresponds a positive number $\beta_{1}$ such that
$P_{a}(\beta^{-1}X)\leq 1-\frac{1}{a}$ $(a>1)$ whenever $\beta>\beta_{1}$. This will be the key operation for proving Theorem \ref{the2}.
\begin{lemma}\label{lem3}
Let $X=U\mathrm{Diag}(\sigma(X))V^{\top}$ be the SVD of matrix $X$, and $rank(X)=r$. Then there exists
\begin{equation}\label{equ18}
\beta_{1} = \displaystyle\frac{a(ar-a+1)\sigma_{1}(X)}{a-1} \ \ (a>1)
\end{equation}
such that, for any $\beta\geq \beta_{1}$,
\begin{equation}\label{equ19}
P_{a}(\beta^{-1}X)\leq 1-\frac{1}{a} \ \ (a>1).
\end{equation}
\end{lemma}
\textbf{proof.} Since the non-convex fraction function $\rho_{a}(t)$ is increasing in $t\in[0,+\infty)$, we have
\begin{equation}\label{equ20}
\begin{array}{llll}
P_{a}(\beta^{-1}X)&=&\displaystyle\sum_{i=1}^{r}\rho_{a}(\sigma_{i}(\beta^{-1}X))\\
&\leq& r\rho_{a}(\sigma_{1}(\beta^{-1}X))\\
&=&\displaystyle\frac{ar\sigma_{1}(X)}{a\sigma_{1}(X)+\beta}.
\end{array}
\end{equation}
In order to get $P_{a}(\beta^{-1}X)\leq \displaystyle 1-\frac{1}{a}$, it suffices to impose
\begin{equation}\label{equ21}
\frac{ar\sigma_{1}(X)}{a\sigma_{1}(X)+\beta}\leq 1-\frac{1}{a},
\end{equation}
equivalently,
$$\beta \geq \frac{a(ar-a+1)\sigma_{1}(X)}{a-1}.$$
This completes the proof. $\hfill{} \Box$

We now proceed to a proof of Theorem \ref{the2}.\\
\textbf{proof.} [\textbf{of Theorem \ref{the2}}]
For each $j\in I_{i}$, combing the definition of $P_{a}$ and Lemma \ref{lem3}, we have
$$\rho_{a}\Big(\sigma_{j}(\gamma^{-1}R_{i})\Big)\leq P(\gamma^{-1}R_{i})\leq1-\frac{1}{a}.$$
Also since
$$\frac{a\sigma_{j}(\gamma^{-1}R_{i})}{a\sigma_{j}(\gamma^{-1}R_{i})+1}\leq 1-\frac{1}{a} \ \  \Leftrightarrow \ \ \sigma_{j}(\gamma^{-1}R_{i})\leq 1-\frac{1}{a}$$
we can get that
$$\sigma_{j}(\gamma^{-1}R_{i})\leq\rho_{a}\Big(\sigma_{j}(\gamma^{-1}R_{i})\Big),\ \ \forall j\in I_{i}.$$
According to the facts that the non-convex fraction function $\rho_{a}(t)$ is increasing for $t>0$, and
$\sigma_{j}(R_{i})\leq\sigma_{k}(R_{i-1})$ for each $j\in I_{i}$ and $k\in I_{i-1}$, $i\geq 2$, we have
$$\sigma_{j}(\gamma^{-1}R_{i})\leq \rho_{a}\Big(\sigma_{j}(\gamma^{-1}R_{i})\Big)\leq\frac{P(\gamma^{-1}R_{i-1})}{K}.$$
It follows that
$$\|\gamma^{-1}R_{i}\|_{F}\leq\frac{P(\gamma^{-1}R_{i-1})}{\sqrt{K}}$$
and
$$\sum_{i\geq2}\|\gamma^{-1}R_{i}\|_{F}\leq\sum_{i\geq2}\frac{P(\gamma^{-1}R_{i-1})}{\sqrt{K}}.$$
Combined with Corollary \ref{cor2}, we immediately get the second part of inequalities (\ref{equ17}).
This completes the proof. $\hfill{} \Box$

\section{The equivalence between the problem (TrAMRM) and (AMRM)} \label{section3}
In this section, a sufficient condition on equivalence of the problem (TrAMRM) and (AMRM) is demonstrated, we proved that
the optimal solution to problem (TrAMRM) also solves (AMRM) if some specific conditions are satisfied.

\begin{definition}\label{de3} {\rm(see \cite{Recht10})}
Let $\mathcal{A}: \mathbb{R}^{m\times n}\mapsto \mathbb{R}^{d} $ be a linear map. For every integer $r$ with
$1\leq r\leq m$, define the $r$-restricted isometry constant to be the smallest number $\delta_{r}(\mathcal{A})$
such that
\begin{equation}\label{equ22}
(1-\delta_{r}(\mathcal{A}))\|X\|_{F}^{2}\leq\|\mathcal{A}(X)\|_{2}^{2}\leq(1+\delta_{r}(\mathcal{A}))\|X\|_{F}^{2}
\end{equation}
holds for all matrix $X\in \mathbb{R}^{m\times n}$ of rank at most $r$.
\end{definition}

Based on Definition \ref{de3}, we shall demonstrate that the optimal solution of the problem (TrAMRM) equivalences to the problem (AMRM).
\begin{theorem}\label{the3}
Let $X^{\ast}$ and $X_{0}$ be the optimal solutions to the problem (TrAMRM) and (AMRM), respectively. If there is a number $K>2T$, such that
\begin{equation}\label{equ23}
\frac{K}{2T}\Big(1-\delta_{2T+K}(\mathcal{A})\Big)-\Big(1+\delta_{K}(\mathcal{A})\Big)>0,
\end{equation}
then there exists $a^{\ast}>1$ (depends on $\delta_{K}(\mathcal{A})$ and $\delta_{2T+K}(\mathcal{A})$), such that for any $1<a<a^{\ast}$, $X^{\ast}=X_{0}$,
where $\mathrm{rank}(X_{0})=T$.
\end{theorem}
\textbf{proof.}
Define the function
$$f(a)=\frac{1}{a^{2}}\frac{K}{2T}\Big(1-\delta_{2T+K}(\mathcal{A})\Big)-1-\delta_{K}(\mathcal{A}) \ \ (a>0).$$
Clearly, the function $f$ is continuous and decreasing in $a\in(0,+\infty)$. Notice that at $a=1$,
\begin{equation}\label{equ24}
f(1)=\frac{K}{2T}\Big(1-\delta_{2T+K}(\mathcal{A})\Big)-1-\delta_{K}(\mathcal{A})>0,
\end{equation}
and as $a\rightarrow+\infty$, $f(a)\rightarrow -1-\delta_{K}(\mathcal{A})<0$. Then, there exists a constant $a^{\ast}>1$ such
that $f(a^{\ast})=0$. It is obvious that the number $a^{\ast}$ depends only on the RIC of linear map $\mathcal{A}$.
Thus, for any $1<a<a^{\ast}$, we have
\begin{equation}\label{equ25}
\frac{1}{a}\sqrt{\frac{1-\delta_{2T+K}}{2T}}-\sqrt{\frac{1+\delta_{K}}{K}}>0.
\end{equation}
Let $R=X^{\ast}-X_{0}$, and in order to show that $X^{\ast}=X_{0}$, it suffices to show that the matrix $R=0$. Partition matrix $R$ as matrices $R_{0}$ and $R_{c}$ which are defined
in Definitions \ref{def1} and \ref{def2}. Since $\mathcal{A}(R)=\mathcal{A}(X^{\ast}-X_{0})=0$, we can get that
\begin{equation}\label{equ26}
\begin{array}{llll}
0&=&\|\mathcal{A}(\gamma^{-1}R)\|_{2}\\
&=&\|\mathcal{A}(\gamma^{-1}R_{0}+\gamma^{-1}R_{c})\|_{2}\\
&=&\|\mathcal{A}(\gamma^{-1}R_{0}+\gamma^{-1}R_{1})+\displaystyle\sum_{i\geq2}\mathcal{A}(\gamma^{-1}R_{i})\|_{2}\\
&\geq&\|\mathcal{A}(\gamma^{-1}R_{0}+\gamma^{-1}R_{1})\|_{2}-\displaystyle\|\sum_{i\geq2}\mathcal{A}(\gamma^{-1}R_{i})\|_{2}\\
&\geq&\|\mathcal{A}(\gamma^{-1}R_{0}+\gamma^{-1}R_{1})\|_{2}-\displaystyle\sum_{i\geq2}\|\mathcal{A}(\gamma^{-1}R_{i})\|_{2}\\
&\geq&\sqrt{1-\delta_{2T+K}(\mathcal{A})}\|\gamma^{-1}R_{0}+\gamma^{-1}R_{1}\|_{F}-\sqrt{1+\delta_{K}(\mathcal{A})}\displaystyle\sum_{i\geq2}\|\gamma^{-1}R_{i}\|_{F}.
\end{array}
\end{equation}
Plus inequalities (\ref{equ16}) and (\ref{equ17}) into inequality (\ref{equ26}), we can get that
\begin{equation}\label{equ27}
\begin{array}{llll}
0&\geq&\displaystyle\sqrt{1-\delta_{2T+K}(\mathcal{A})}\frac{1}{a\sqrt{2T}}P_{a}(\gamma^{-1}R_{0})-\sqrt{1+\delta_{K}(\mathcal{A})}\frac{1}{\sqrt{K}}P_{a}(\gamma^{-1}R_{0})\\
&=&\displaystyle\Big(\frac{1}{a}\sqrt{\frac{1-\delta_{2T+K}(\mathcal{A})}{2T}}-\sqrt{\frac{1+\delta_{K}(\mathcal{A})}{K}}\big)P_{a}(\gamma^{-1}R_{0}) .
\end{array}
\end{equation}
Moreover, following the inequality (\ref{equ25}), the factor
$$\frac{1}{a}\sqrt{\frac{1-\delta_{2T+K}(\mathcal{A})}{2T}}-\sqrt{\frac{1+\delta_{K}(\mathcal{A})}{K}}$$
is strictly positive for any $1<a<a^{\ast}$, and thus $P_{a}(\gamma^{-1}R_{0})=0$, which implies that $R_{0}=O$. Combined with Corollary \ref{cor2}, $R_{c}=O$.
Therefore, $X^{\ast}=X_{0}$. This completes the proof.$\hfill{} \Box$

\begin{corollary}\label{cor3}
Suppose that the positive integer $T\geq 1$ is such that $\delta_{5T}(\mathcal{A})<\frac{3-2a^{2}}{3+2a^{2}}$ for any $a>1$, then $X^{\ast}=X_{0}$.
\end{corollary}
\textbf{proof.}
By Definition 1, $\delta_{r_{1}}(\mathcal{A})\leq\delta_{r_{2}}(\mathcal{A})$ for $r_{1}\leq r_{2}$. Let $K=3T$, notice that the inequality (\ref{equ25}) holds
when $\frac{3}{2a^{2}}(1-\delta_{5T}(\mathcal{A}))>1+\delta_{3T}(\mathcal{A})$. Since $\delta_{3T}(\mathcal{A})\leq\delta_{5T}(\mathcal{A})$, we immediately
get that $X^{\ast}=X_{0}$ if $\delta_{5T}(\mathcal{A})<\frac{3-2a^{2}}{3+2a^{2}}$. This completes the proof. $\hfill{} \Box$

Theorem \ref{the3} or Corollary \ref{cor3} demonstrated that the optimal solution to the problem (AMRM) can be exactly obtained by solving problem (TrAMRM) if some
specific conditions satisfied. Moreover, we also proved that the optimal solution to the problem (TrAMRM) can be approximately obtained by solving problem (RTrAMRM) for
some proper smaller $\lambda>0$.
\begin{theorem}\label{the4}
Let $\{\lambda_{\tilde{n}}\}$ be a decreasing sequence of
positive numbers with $\lambda_{\tilde{n}}\rightarrow 0$, and $X_{\lambda_{\tilde{n}}}$ be the optimal solution of the problem (RTrAMRM) with
$\lambda=\lambda_{\tilde{n}}$. If the problem (TrAMRM) is feasible, then the sequence $\{X_{\lambda_{\tilde{n}}}\}$ is bounded and any of
its accumulation points is the optimal solution of the problem (TrAMRM).
\end{theorem}
\textbf{proof.} By
$$\lambda_{\tilde{n}} P_{a}(X)\leq\|\mathcal{A}(X)-b\|_{2}^{2}+\lambda_{\tilde{n}} P_{a}(X),$$
we can get that the objective function in the problem (RTrAMRM) with $\lambda=\lambda_{\tilde{n}}$
is bounded from below and is coercive, i.e.,
$$\|\mathcal{A}(X)-b\|_{2}^{2}+\lambda_{\tilde{n}} P_{a}(X)\rightarrow+\infty\ \ \mathrm{as}\ \ \|X\|_{F}\rightarrow+\infty,$$
and hence the set of optimal solution of the problem (RTrAMRM) with $\lambda=\lambda_{\tilde{n}}$ is nonempty and bounded.

By assumption, we suppose that the problem (TrAMRM) is feasible and $\bar{X}$ is any feasible point, then
$$\mathcal{A}(\bar{X})=b.$$
Since $\{X_{\lambda_{\tilde{n}}}\}$ is the optimal solution of the problem (RTrAMRM) with $\lambda=\lambda_{\tilde{n}}$, we have
\begin{equation}\label{equ28}
\begin{array}{llll}
\lambda_{\tilde{n}}P_{a}(X_{\lambda_{\tilde{n}}})&\leq&\|\mathcal{A}(X_{\lambda_{\tilde{n}}})-b\|_{2}^{2}+\lambda_{\tilde{n}}P_{a}(X_{\lambda_{\tilde{n}}})\\
&\leq&\|\mathcal{A}(\bar{X})-b\|_{2}^{2}+\lambda_{\tilde{n}} P_{a}(\bar{X})\\
&=&\lambda_{\tilde{n}} P_{a}(\bar{X}).
\end{array}
\end{equation}
Hence, the sequence $\{P_{a}(X_{\lambda_{\tilde{n}}})\}_{\tilde{n}\in N^{+}}$ is bounded, and the sequence $\{X_{\lambda_{\tilde{n}}}\}$ has at least one
accumulation point. In addition, by inequality (\ref{equ28}), we can get that
$$\|\mathcal{A}(X_{\lambda_{\tilde{n}}})-b\|_{2}^{2}\leq\lambda_{\tilde{n}} P_{a}(\bar{X})\ \ \mathrm{for}\ \mathrm{any}\ \ \lambda_{\tilde{n}}\rightarrow 0.$$
If we set $X^{\ast}$ be any accumulation point of the sequence $\{X_{\lambda_{\tilde{n}}}\}$, we can derive that
$$\mathcal{A}(X^{\ast})=b.$$
That is, $X^{\ast}$ is a feasible point of the problem (TrAMRM). Combined with $P_{a}(X^{\ast})\leq P_{a}(\bar{X})$ and the arbitrariness of $\bar{X}$,
we can get that $X^{\ast}$ is the optimal solution of the problem (TrAMRM). This completes the proof. $\hfill{} \Box$

\section{Algorithm for solving the problem (RTrAMRM)} \label{section4}
In this section, the DC (Difference of Convex functions) algorithm is utilized to solve the non-convex problem (RTrAMRM). For the sake of simplicity, we call it as RTrDC algorithm.

\subsection{{DC programming and DC algorithm} \label{condition-sec}}

\begin{definition}\label{de2} {\rm(DC functions \cite{Thi18,Thi19})}
Let $\mathcal{C}$ be a convex subset of $\mathbb{R}^{l}$. A real-valued function $f:\mathcal{C}\mapsto \mathbb{R}$ is called DC (Difference of Convex functions)
on $\mathcal{C}$, if there exist two convex functions $g,h: \mathcal{C}\mapsto \mathbb{R}$ such that $f$ can be expressed in the form
\begin{equation}\label{equ29}
f(x)=g(x)-h(x).
\end{equation}
If $\mathcal{C}=\mathbb{R}^{l}$, then $f$ is simply called a DC function. Each representation of the form (\ref{equ29}) is said to be a DC decomposition of $f$.
\end{definition}

Generally speaking, the DC programming is an optimization problem of the form
$$\alpha=\inf_{x\in \mathbb{R}^{l}}\{f(x)=g(x)-h(x)\}$$
where $g$, $h$ are lower semi-continuous proper convex functions on $\mathbb{R}^{l}$. The main ideal of DC algorithm is to replace in the DC programming, at the
current point $x^{k}$ of iteration $k$, the second component $h$ with its affine minimization defined by
\begin{equation}\label{equ30}
h_{k}(x)=h(x^{k})+\langle x-x^{k}, y^{k}\rangle, \ \ \ y^{k}\in\partial h(x^{k})
\end{equation}
to give birth to the convex programming of the form
\begin{equation}\label{equ31}
\inf_{x\in \mathbb{R}^{l}}\{g(x)-h_{k}(x)\}\Leftrightarrow \inf_{x\in \mathbb{R}^{l}}\{g(x)-\langle x, y^{k}\rangle\}
\end{equation}
whose optimal solution is taken as $x^{k+1}$.

\begin{algorithm}[h!]
\caption{: DC algorithm}
\label{alg:1}
\begin{algorithmic}
\STATE {\textbf{Initialize}: Let $x^{0}\in \mathbb{R}^{l}$ be an initial guess;}
\STATE {$k=0$;}
\STATE {\textbf{Repeat}}
\STATE \ \ \ \  {$y^{k}\in \partial h(x^{k})$}
\STATE \ \ \ \  {$x^{k+1}\in\arg \displaystyle\min_{x\in \mathbb{R}^{l}}\{g(x)-\langle x,y^{k}\rangle\}$}
\STATE {$k\rightarrow k+1$}
\STATE{\textbf{Until} convergence of $\{x^{k}\}$.}
\end{algorithmic}
\end{algorithm}

\subsection{DC algorithm for solving the problem (RTrAMRM)}
Let
\begin{equation}\label{equ32}
\mathcal{T}_{\lambda}(X)=\|\mathcal{A}(X)-b\|_{2}^{2}+\lambda P_{a}(X).
\end{equation}
It is to see clear that the function $\mathcal{T}_{\lambda}(X)$ is a DC function of the form
$$\mathcal{T}_{\lambda}(X)=g(X)-h(X),$$
where
$$g(X)=\|\mathcal{A}(X)-b\|_{2}^{2}+\lambda \|X\|_{\ast}$$
and
$$h(X)=\lambda \|X\|_{\ast}-\lambda P(X)$$
are all convex functions. Hence, the resulting problem of (RTrAMRM) via this approximation can be written as a DC program
\begin{equation}\label{equ33}
\min_{X\in \mathbb{R}^{m\times n}}\Big\{(\|\mathcal{A}(X)-b\|_{2}^{2}+\lambda \|X\|_{\ast})-(\lambda \|X\|_{\ast}-\lambda P_{a}(X))\Big\}.
\end{equation}
Applying DC algorithm on (\ref{equ33}) amounts to computing the two sequences $\{N^{k}\}$ and $\{X^{k}\}$ such that $N^{k}\in\partial (\lambda \|X\|_{\ast}-\lambda P_{a}(X))$
and $X^{k+1}$ is the solution to the following convex problem
\begin{equation}\label{equ34}
X^{k+1}\in\arg \min_{X\in \mathbb{R}^{m\times n}}\Big\{\|\mathcal{A}(X)-b\|_{2}^{2}+\lambda \|X\|_{\ast}-\langle X,N^{k}\rangle\Big\}.
\end{equation}
It is necessary to emphasize that, at each iteration, we need to solve a convex sub-problem (\ref{equ34}).

\begin{algorithm}[h!]
\caption{: RTrDC algorithm}
\label{alg:2}
\begin{algorithmic}
\STATE {\textbf{Input}: $\mathcal{A}: \mathbb{R}^{m\times n}\mapsto \mathbb{R}^{d}$, $b\in \mathbb{R}^{d}$;}
\STATE {\textbf{Initialize}: Given $X^{0}\in \mathbb{R}^{m\times n}$, $a>0$ and $\lambda>0$;}
\STATE {$k=0$;}
\STATE {\textbf{Repeat}}
\STATE \ \ \ \  {$N^{k}\in \partial(\lambda \|X^{k}\|_{\ast}-\lambda P_{a}(X^{k}))$}
\STATE \ \ \ \  {$X^{k+1}\in\displaystyle\arg\min_{X\in \mathbb{R}^{m\times n}}\Big\{\|\mathcal{A}(X)-b\|_{2}^{2}+\lambda \|X\|_{\ast}-\langle X,N^{k}\rangle\Big\}$}
\STATE {$k\rightarrow k+1$}
\STATE{\textbf{Until} convergence of $\{X^{k}\}$.}
\end{algorithmic}
\end{algorithm}

\begin{remark}\label{re1}
Let $X^{k}=U^{k}\mathrm{Diag}(\sigma_{i}(X^{k}))V^{k}$ be the SVD of the matrix $X^{k}$. Then $\partial(\lambda \|X^{k}\|_{\ast}-\lambda P_{a}(X^{k}))=U^{k}\mathrm{Diag}\Big(\lambda-\frac{\lambda b}{(b\sigma_{i}(X^{k})+1)^{2}}\Big)V^{k}$. The detailed proof can be seen in \cite{Lew20}.
\end{remark}

Before continuing our discussion, the definition of the singular value thresholding operator \cite{Cai14} should be prepared, which underlies the closed
form representation of the optimal solution to the problem (\ref{equ34}).

\begin{definition}\label{def5}{\rm(see \cite{Cai14})}
Let $Y=U\Sigma V^{\top}=U\mathrm{Diag}(\sigma_{i}(Y)) V^{\top}$ be the SVD of matrix $Y$, for any $\lambda>0$, suppose that
\begin{equation}\label{equ35}
\mathcal{D}_{\lambda}(Y)\triangleq\arg\min_{X\in \mathbb{R}^{m\times n}}\Big\{\|X-Y\|_{F}^{2}+\lambda\|X\|_{\ast}\Big\},
\end{equation}
then the soft thresholding operator $\mathcal{D}_{\lambda}$ can be expressed as
\begin{equation}\label{equ36}
\mathcal{D}_{\lambda}(Y)=U\mathcal{D}_{\lambda}(\Sigma)V^{\top}=U\mathrm{Diag}\Big(\{\sigma_{i}(Y)-\frac{\lambda}{2}\}_{+}\Big)V^{\top}
\end{equation}
where $t_{+}$ is the positive part of $t$, and $t_{+}=\max(0,t)$.
\end{definition}

The singular value thresholding operator $\mathcal{D}_{\lambda}$ simply applies the soft thresholding operator \cite{Dau21} defined on vector to the
singular values of a matrix, and effectively shrinks them towards zero. In particular, it needs to be emphasized that the soft thresholding operator has
been actively studied in different fields such as signal processing \cite{Dau21,Yu22}, statistics \cite{Tib23}, portfolio section \cite{Jos24} and
visual tracking \cite{Lan25,Lan26,Lan27,Lan29,Lan30}.

Nextly, we will show that the optimal solution to the problem (\ref{equ34}) can be expressed a thresholding operation.

Let
$$\mathcal{L}_{1}(X)=\|\mathcal{A}(X)-b\|_{2}^{2}+\lambda\|X\|_{\ast}-\langle X,N^{k}\rangle$$
and its surrogate function
\begin{eqnarray*}
\mathcal{L}_{2}(X,Z,\mu)&=&\mu\mathcal{L}_{1}(X)-\mu\|\mathcal{A}(X)-\mathcal{A}(Z)\|_{2}^{2}+\|X-Z\|_{F}^{2}
\end{eqnarray*}
where $Z\in \mathbb{R}^{m\times n}$ is an additional variable. Clearly $\mathcal{L}_{2}(X,X,\mu)=\mu\mathcal{L}_{1}(X)$.

\begin{theorem}\label{the5}
For any fixed $\lambda>0$, $\mu>0$ and $Z\in \mathbb{R}^{m\times n}$, $\displaystyle\min_{X\in \mathbb{R}^{m\times n}}\mathcal{L}_{2}(X,Z,\mu)$ equivalents to
\begin{equation}\label{equ37}
\min_{X\in \mathbb{R}^{m\times n}}\Big\{\|X-B_{\mu}(Z)\|_{F}^{2}+\lambda\mu \|X\|_{\ast}\Big\},
\end{equation}
where $B_{\mu}(Z)=Z+\mu \mathcal{A}^{\ast}(b-\mathcal{A}(Z))+\frac{1}{2}\mu N^{k}$.
\end{theorem}
\textbf{proof.} By the definition of $\mathcal{L}_{2}(X,Z,\mu)$, we have
\begin{eqnarray*}
\mathcal{L}_{2}(X,Z,\mu)&=&\mu\|\mathcal{A}(X)-b\|_{2}^{2}+\lambda\mu\|X\|_{\ast}-\mu\langle X,N^{k}\rangle-\mu\|\mathcal{A}(X)-\mathcal{A}(Z)\|_{2}^{2}\\
&&+\|X-Z\|_{F}^{2}\\
&=&\|X-(Z+\mu \mathcal{A}^{\ast}(b-\mathcal{A}(Z))+\frac{1}{2}\mu N^{k})\|_{F}^{2}+\lambda\mu\|X\|_{\ast}+\|Z\|_{F}^{2}\\
&&-\|Z+\mu \mathcal{A}^{\ast}(b-\mathcal{A}(Z))+\frac{1}{2}\mu N^{k}\|_{F}^{2}+\mu\|b\|_{2}^{2}-\mu\|\mathcal{A}(Z)\|_{2}^{2}\\
&=&\|X-B_{\mu}(Z)\|_{F}^{2}+\lambda\mu\|X\|_{\ast}+\|Z\|_{F}^{2}-\|B_{\mu}(Z)\|_{F}^{2}+\mu\|b\|_{2}^{2}\\
&&-\mu\|\mathcal{A}(Z)\|_{2}^{2}.
\end{eqnarray*}
This completes the proof. $\hfill{} \Box$

\begin{theorem}\label{the6}
For fixed positive parameters $\lambda>0$ and $0<\mu<\frac{1}{\|\mathcal{A}\|_{2}^{2}}$. If the matrix $X^{\ast}$ is the optimal solution
to $\displaystyle\min_{X\in \mathbb{R}^{m\times n}}\mathcal{L}_{1}(X)$, then $X^{\ast}$ is also the optimal solution to
$\displaystyle\min_{X\in \mathbb{R}^{m\times n}}\mathcal{L}_{2}(X,X^{\ast},\mu)$, that is
$$\mathcal{L}_{2}(X^{\ast},X^{\ast},\mu)\leq \mathcal{L}_{2}(X,X^{\ast},\mu)$$
for any $X\in \mathbb{R}^{m\times n}$.
\end{theorem}
\textbf{proof.} By the definition of $\mathcal{L}_{2}(X,Z,\mu)$, we have
\begin{eqnarray*}
\mathcal{L}_{2}(X,X^{\ast},\mu)&=&\mu\|\mathcal{A}(X)-b\|_{2}^{2}+\lambda\mu\|X\|_{\ast}-\mu\langle X,N^{k}\rangle\\
&&-\mu\|\mathcal{A}(X)-\mathcal{A}(X^{\ast})\|_{2}^{2}+\|X-X^{\ast}\|_{F}^{2}\\
&\geq&\mu\|\mathcal{A}(X)-b\|_{2}^{2}+\lambda\mu\|X\|_{\ast}-\mu\langle X,N^{k}\rangle\\
&=&\mu \mathcal{L}_{1}(X)\\
&\geq&\mu \mathcal{L}_{1}(X^{\ast})\\
&=&\mathcal{L}_{2}(X^{\ast},X^{\ast},\mu)
\end{eqnarray*}
where the first inequality holds by the fact that
\begin{eqnarray*}
\|\mathcal{A}(X)-\mathcal{A}(X^{\ast})\|_{2}^{2}&=&\|Avec(X)-Avec(X^{\ast})\|_{2}^{2}\\
&\leq&\|\mathcal{A}\|_{2}^{2}\|X-X^{\ast}\|_{F}^{2}.
\end{eqnarray*}
This completes the proof.$\hfill{} \Box$

Theorem \ref{the6} demonstrated that the matrix $X^{\ast}$ is the global optimal solution to $\displaystyle\min_{X\in \mathbb{R}^{m\times n}}\mathcal{L}_{2}(X,X^{\ast},\mu)$
if and only if the matrix $X^{\ast}$ is the global optimal solution to $\displaystyle\min_{X\in \mathbb{R}^{m\times n}}\mathcal{L}_{1}(X)$. Combing with Theorem \ref{the5},
we can get that the optimal solution to $\displaystyle\min_{X\in \mathbb{R}^{m\times n}}\mathcal{L}_{2}(X,X^{\ast},\mu)$ could be obtained by
solving the following problem:
\begin{equation}\label{equ38}
\min_{X\in \mathbb{R}^{m\times n}}\Big\{\|X-B_{\mu}(X^{\ast})\|_{F}^{2}+\lambda\mu\|X\|_{\ast}\Big\}.
\end{equation}
where $B_{\mu}(X^{\ast})=X^{\ast}+\mu \mathcal{A}^{\ast}(b-\mathcal{A}(X^{\ast}))+\frac{1}{2}\mu N^{k}$. Moreover, by Definition \ref{def5}, the optimal solution to the minimization
 problem (\ref{equ38}) could be deduced to the following form
\begin{equation}\label{equ39}
X^{\ast}=\mathcal{D}_{\lambda\mu}(B_{\mu}(X^{\ast}))=U^{\ast}\mathcal{D}_{\lambda\mu}(\Sigma_{B}^{\ast})(V^{\ast})^{\top}
\end{equation}
where $B_{\mu}(X^{\ast})=U^{\ast}\Sigma^{\ast}(V^{\ast})^{\top}=U^{\ast}\mathrm{Diag}(\sigma_{i}(B_{\mu}(X^{\ast})))(V^{\ast})^{\top}$
is the SVD of the matrix $B_{\mu}(X^{\ast})$, and the operator $\mathcal{D}_{\lambda\mu}$ is obtained by replacing $\lambda$ with $\lambda\mu$
in $\mathcal{D}_{\lambda}$.

With the thresholding representation (\ref{equ39}), the procedure of the thresholding algorithm for solving the sub-problem (\ref{equ34})
can be naturally defined as
\begin{equation}\label{equ40}
\begin{array}{llll}
X^{s+1}&=&\mathcal{D}_{\lambda\mu}(B_{\mu}(X^{s}))\\
&=&U^{s}\mathrm{Diag}\Big(\{\sigma_{i}(B_{\mu}(X^{s}))-\frac{\lambda\mu}{2}\}_{+}\Big)(V^{s})^{\top}
\end{array}
\end{equation}
until a stopping criterion is reached.

It is necessary to emphasize that the quantity of the solution of a regularization problem depends seriously on the setting of the regularization parameter $\lambda>0$, and the
selection of proper regularization parameter is a very hard problem. In iteration (\ref{equ40}), the cross-validation method is accepted to choose the proper regularization
parameter $\lambda>0$. To make it clear, we suppose that the matrix $X^{\ast}$ of rank $r$ is the optimal solution to the problem (\ref{equ34}), and the singular values of matrix
$B_{\mu}(X^{\ast})$ are denoted as
$$\sigma_{1}(B_{\mu}(X^{\ast}))\geq\sigma_{2}(B_{\mu}(X^{\ast}))\geq\cdots\geq\sigma_{m}(B_{\mu}(X^{\ast})).$$
By (\ref{equ36}), it then follows that
$$\sigma_{i}(B_{\mu}(X^{\ast}))>\frac{\lambda^{\ast}\mu}{2}\ \Leftrightarrow \ i\in\{1,2,\cdots,r\},$$
$$\sigma_{i}(B_{\mu}(X^{\ast}))\leq\frac{\lambda^{\ast}\mu}{2}\ \Leftrightarrow \ i\in\{r+1,r+2,\cdots,m\},$$
which implies
$$\frac{2\sigma_{r+1}(B_{\mu}(X^{\ast}))}{\mu}\leq\lambda^{\ast}<\frac{2\sigma_{r}(B_{\mu}(X^{\ast}))}{\mu},$$
namely
\begin{equation}\label{equ41}
\lambda^{\ast}\in\bigg[\frac{2\sigma_{r+1}(B_{\mu}(X^{\ast}))}{\mu}, \frac{2\sigma_{r}(B_{\mu}(X^{\ast}))}{\mu}\bigg).
\end{equation}
We can then take
\begin{equation}\label{equ42}
\lambda^{\ast}=\frac{2(1-\theta)\sigma_{r+1}(B_{\mu}(X^{\ast}))}{\mu}+\frac{2\theta\sigma_{r}(B_{\mu}(X^{\ast}))}{\mu}.
\end{equation}
with any $\theta\in[0,1)$. Taking $\theta=0$, this leads to a most reliable choice of $\lambda^{\ast}$ specified by
\begin{equation}\label{equ43}
\lambda^{\ast}=\frac{2\sigma_{r+1}(B_{\mu}(X^{\ast}))}{\mu}.
\end{equation}
In practice, we approximate $B_{\mu}(X^{\ast})$ by $B_{\mu}(X^{s})$ in (\ref{equ43}), and the regularization parameter $\lambda$ could be selected as
\begin{equation}\label{equ44}
\lambda_{s}^{\ast}=\frac{2\sigma_{r+1}(B_{\mu}(X^{s}))}{\mu}
\end{equation}
in applications. When so doing, our algorithm will be adaptive and free from the choice of regularization parameter.

\section{Numerical experiments} \label{section5}
\begin{figure}[h!]
  \centering
  \begin{minipage}[t]{0.4\linewidth}
  \centering
  \includegraphics[width=1\textwidth]{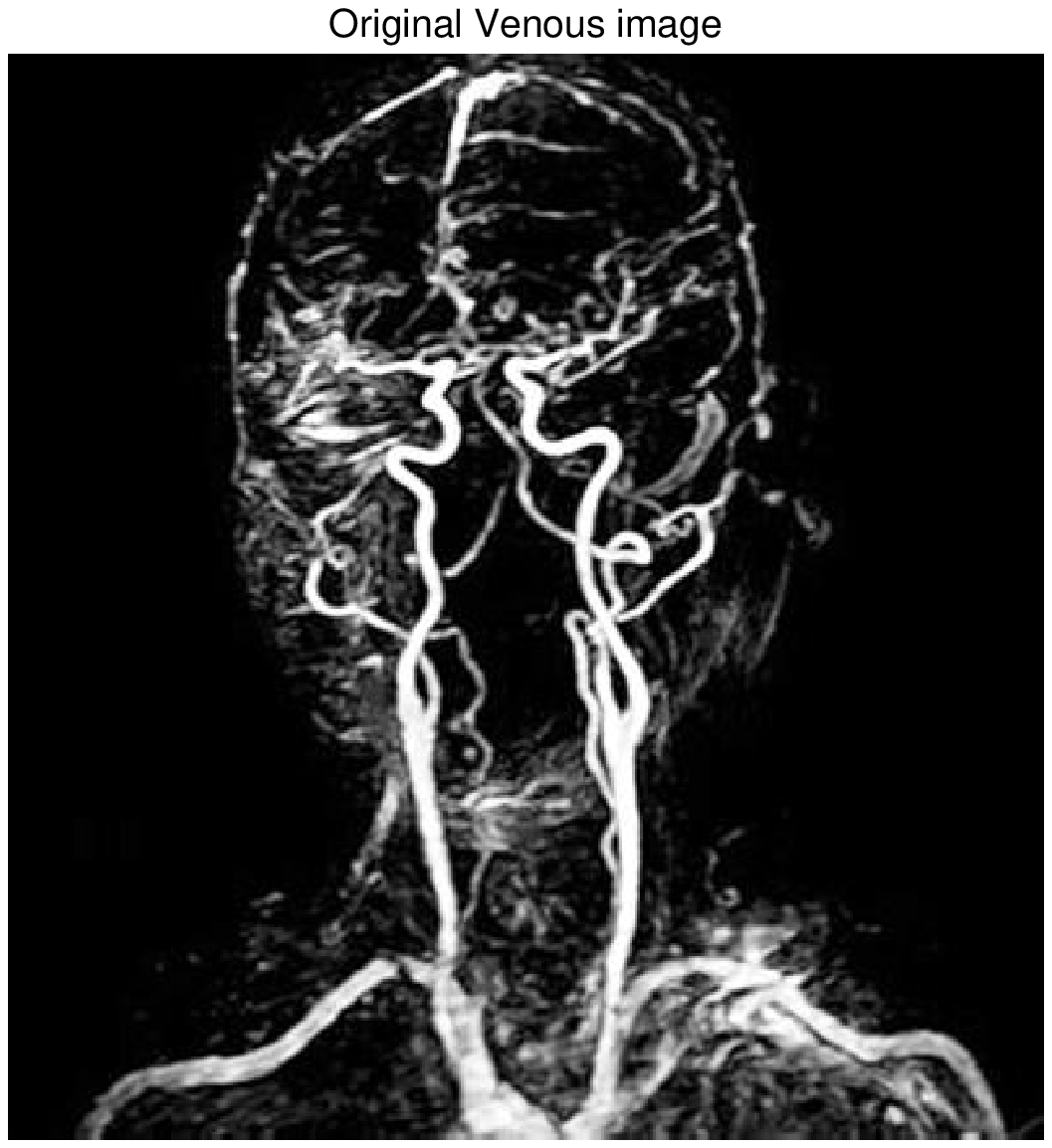}
  \end{minipage}
  \begin{minipage}[t]{0.4\linewidth}
  \centering
  \includegraphics[width=1\textwidth]{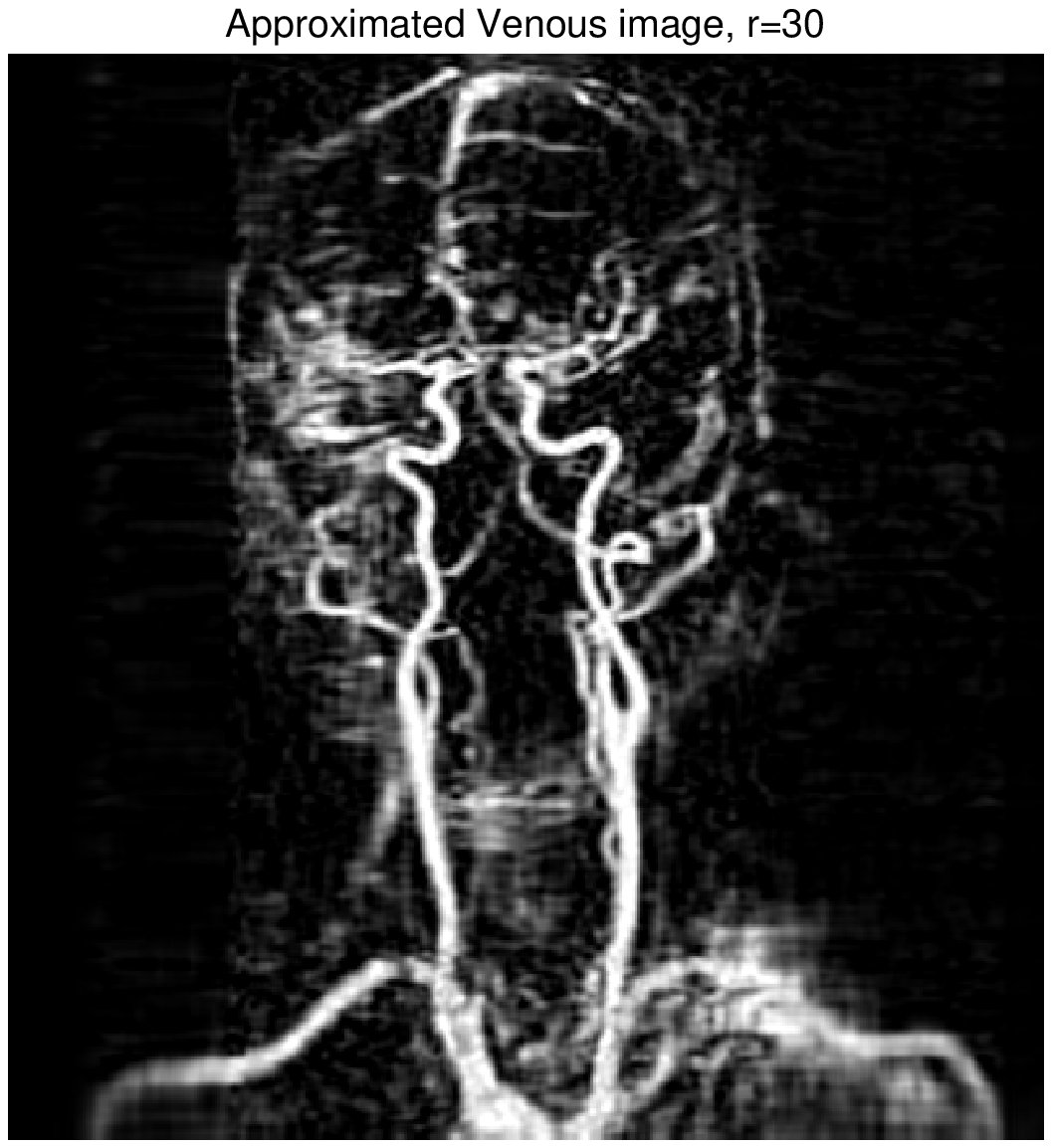}
  \end{minipage}
  \caption{Original $419\times 400$ Venous image and its approximation image with rank $r=30$.} \label{figure1}
\end{figure}

\begin{figure}[h!]
  \centering
  \begin{minipage}[t]{0.4\linewidth}
  \centering
  \includegraphics[width=1\textwidth]{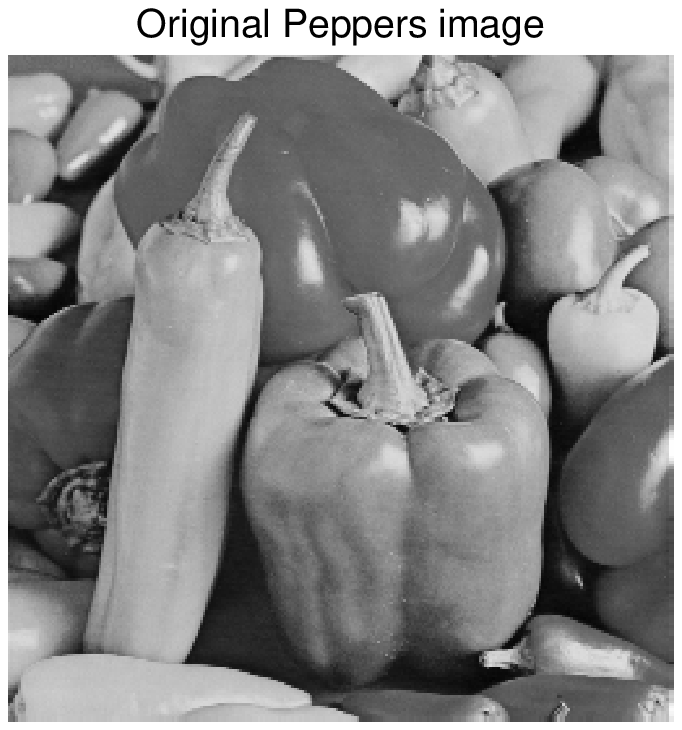}
  \end{minipage}
  \begin{minipage}[t]{0.4\linewidth}
  \centering
  \includegraphics[width=1\textwidth]{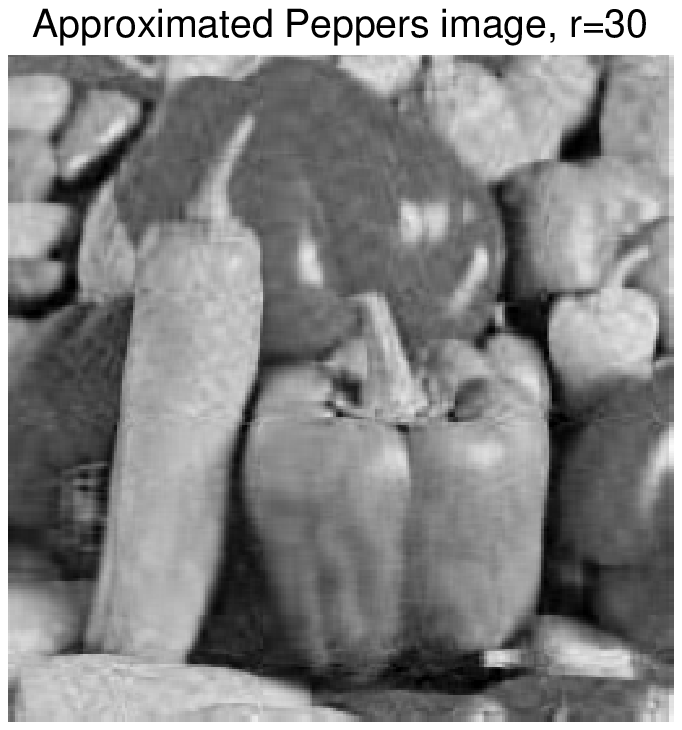}
  \end{minipage}
  \caption{Original $256\times 256$ Peppers image and its approximation image with rank $r=30$.} \label{figure2}
\end{figure}

In this section, we present numerical results of the RTrDC algorithm and compare them with some state-of-art methods (singular value thresholding (SVT)
algorithm \cite{Cai14} and singular value projection (SVP) algorithm \cite{Rag31}) on two image inpainting problems. The three algorithms are tested on
two gray-scale images: $419\times400$ Venous and $256\times 256$ Peppers. We use the SVD to obtain their approximated low rank images with rank $r=30$.
The original images, and their low-rank images are displayed in Figure \ref{figure1} and Figure \ref{figure2} respectively.
The set of observed entries
$\Omega$ is sampled uniformly at random among all sets of cardinality $s$. $\mathrm{SR}=s/mn$ denotes the sampling ration. $\mathrm{FR}=s/r(m+n-r)$ denotes
the freedom ration is the ratio between the number of sampled entries and the `true dimensionality' of a $m\times n$ matrix of rank $r$. If $\mathrm{FR}<1$,
it is impossible to recover an original low-rank matrix because there are an infinite number of matrices of rank $r$ with the observed entries \cite{Ma15}.
The stopping criterion is usually as follows
$$\frac{\|X^{k}-X^{k-1}\|_{F}}{\|X^{k}\|_{F}}\leq \mathrm{Tol}$$
where $X^{k}$ and $X^{k-1}$ are numerical results from two continuous iterative steps and $\mathrm{Tol}$ is a given small number. We take $\mathrm{Tol}=10^{-8}$ in our
experiments. In addition, we measure the accuracy of the generated solution $X^{opt}$ of our algorithms by the relative error ($\mathrm{RE}$) defined as follows
$$\mathrm{RE}=\frac{\|X^{opt}-M\|_{F}}{\|M\|_{F}}.$$
In Theorem \ref{the4}, we have proved that, for some proper smaller $\lambda>0$, the optimal solution to the problem (TrAMRM) can be approximately obtained by solving
problem (RTrAMRM). Moreover, in Theorem \ref{the3} we have proved that there exists $a^{\ast}>1$ such that the optimal solution to the problem (TrAMRM) also solves
the problem (AMRM) whenever $1<a<a^{\ast}$. However, the value of $a^{\ast}$ is extremely difficult to evaluate, and it seriously depends on the rank of the optimal
solution of (AMRM). For the sake of simplicity, in these experiments, we set $a=1.2$, which is closes to 1.

\subsection{Image inpainting-noiseless case }

In this section, we consider the noiseless case and take a series of experiments to demonstrate the performance of the RTrDC algorithm on two image inpainting problems.
\begin{figure}[h!]
  \centering
  \begin{minipage}[t]{0.4\linewidth}
  \centering
  \includegraphics[width=1\textwidth]{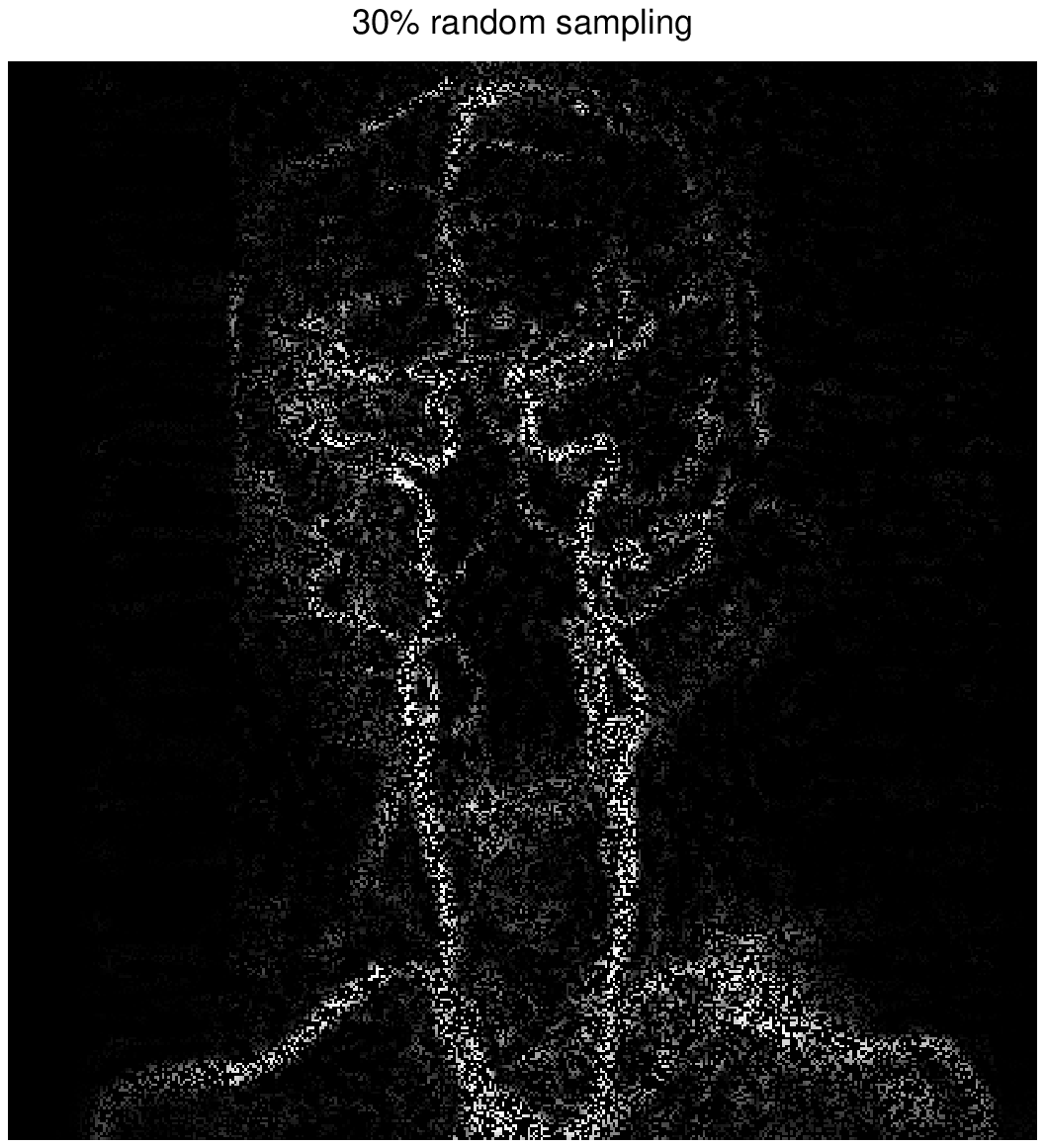}
  \end{minipage}
  \begin{minipage}[t]{0.4\linewidth}
  \centering
  \includegraphics[width=1\textwidth]{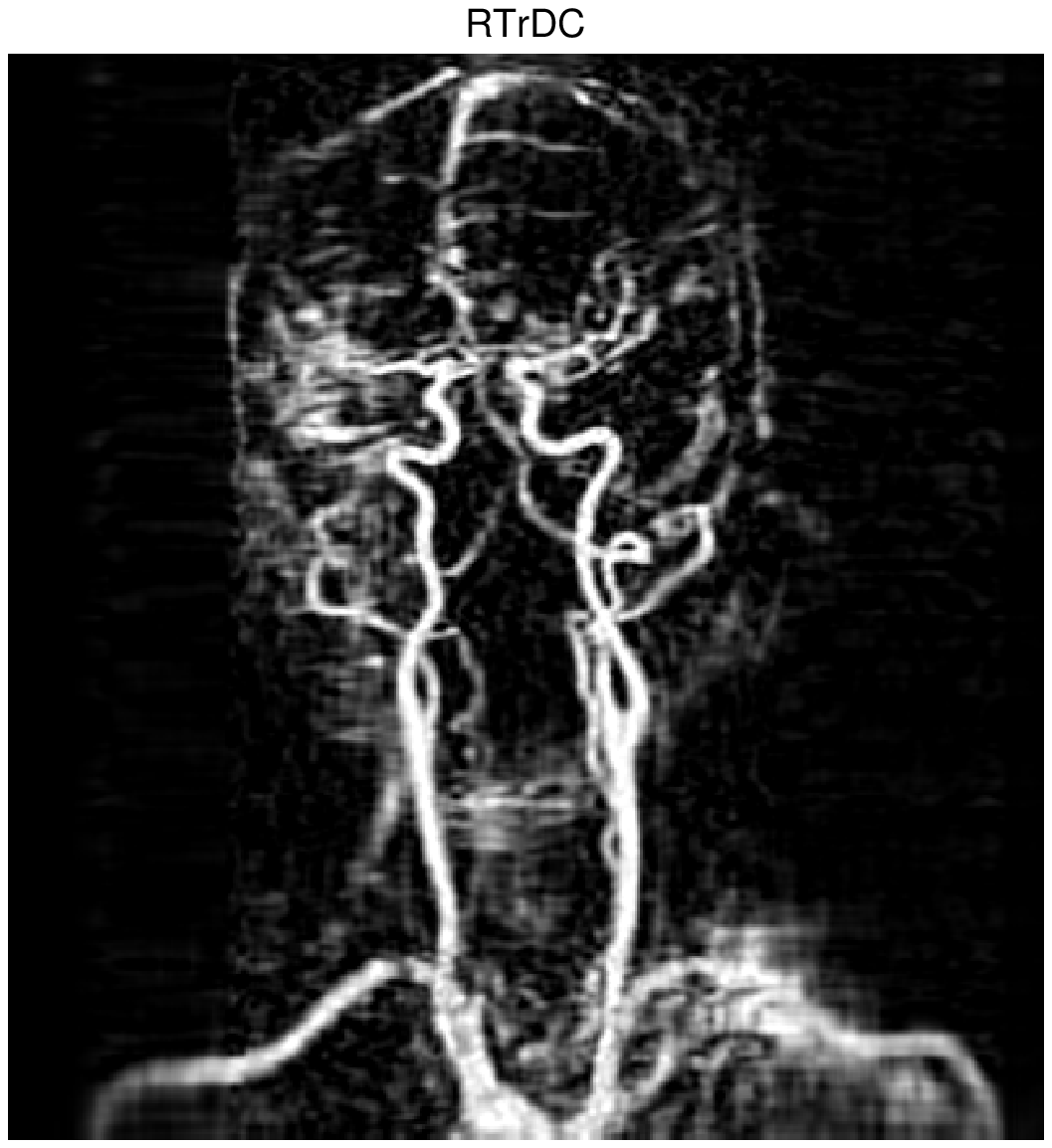}
  \end{minipage}
    \begin{minipage}[t]{0.4\linewidth}
  \centering
  \includegraphics[width=1\textwidth]{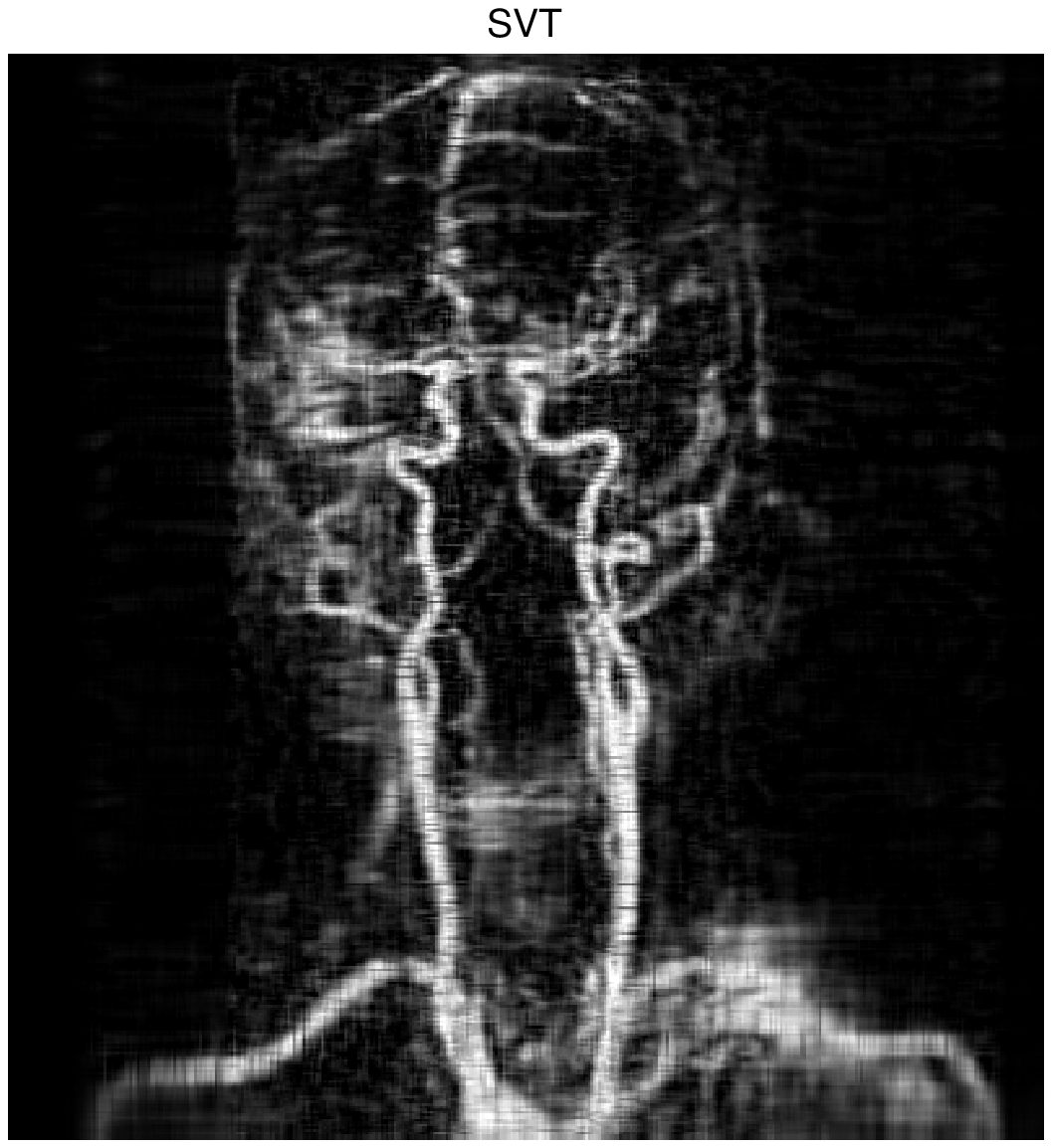}
  \end{minipage}
    \begin{minipage}[t]{0.4\linewidth}
  \centering
  \includegraphics[width=1\textwidth]{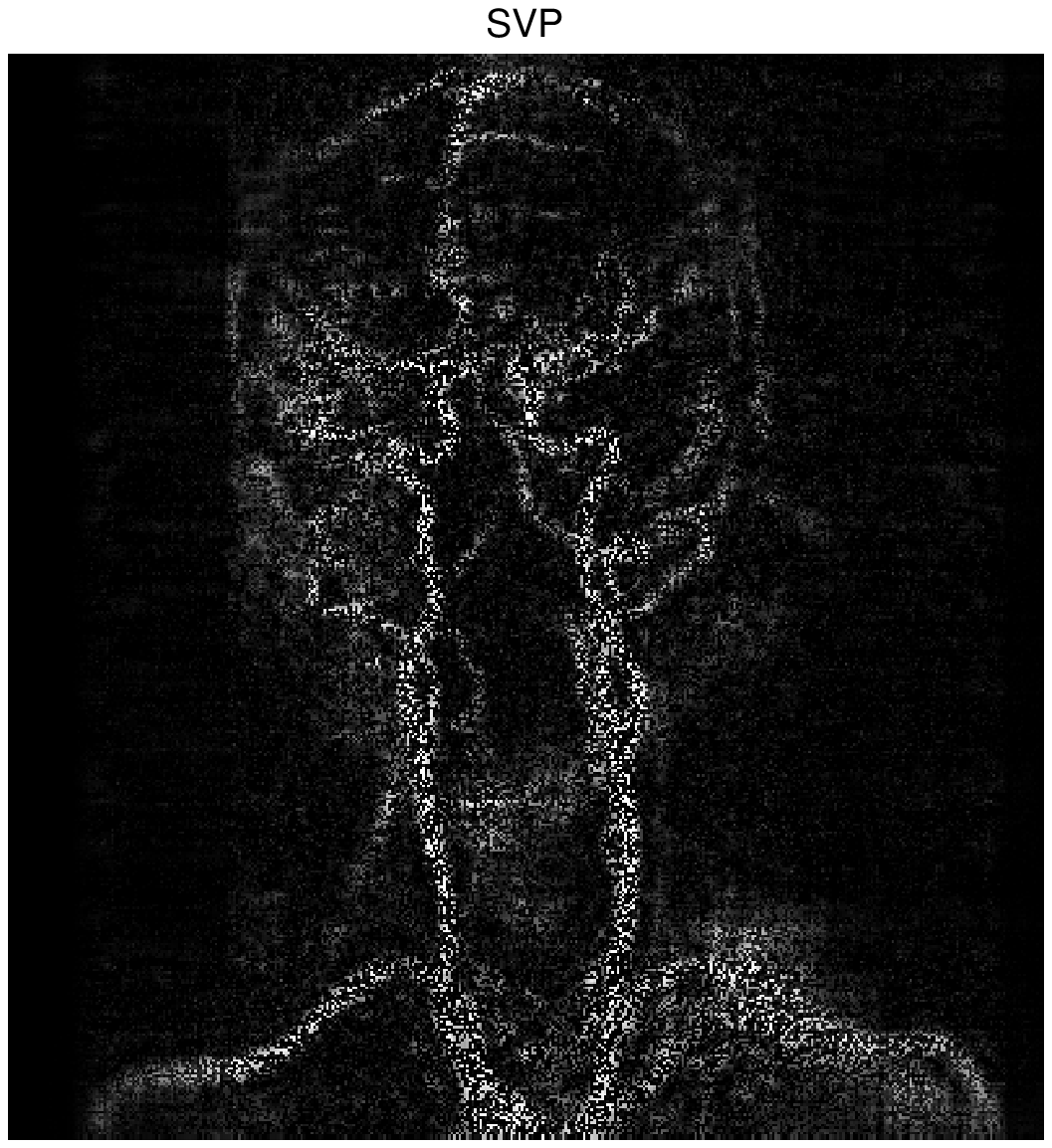}
  \end{minipage}
  \caption{Comparisons of RTrDC algorithm, SVT algorithm and SVP algorithm for recovering the approximated Venous image (noiseless case) with $\mathrm{SR}=0.30$.} \label{figure3}
\end{figure}
\begin{table}[h!]
\footnotesize
\centering
\begin{tabular}{|c||l|l|l|l|}\hline
\multicolumn{4}{|c|}{Venous image, noiseless, $(r=30,\mathrm{SR}=0.40,\mathrm{FR}=2.8323)$}\\\hline
Algorithm & RTrDC algorithm& SVT algorithm& SVP algorithm\\
\hline
RE& 3.52e--06 & 7.29e--02 &7.40e--01\\
\hline
\multicolumn{4}{|c|}{Peppers image, noiseless, $(r=30,\mathrm{SR}=0.40,\mathrm{FR}=1.8129)$}\\\hline
Algorithm & RTrDC algorithm& SVT algorithm& SVP algorithm\\
\hline
RE & 1.81e--05& 4.43e--02& 7.60e--01\\
\hline
\end{tabular}
\caption{\scriptsize Numerical results of RTrDC algorithm, SVT algorithm, SVP algorithm for image inpainting problems (noiseless case), $\mathrm{SR}=0.40$.}\label{table1}
\end{table}

\begin{figure}[h!]
  \centering
  \begin{minipage}[t]{0.4\linewidth}
  \centering
  \includegraphics[width=1\textwidth]{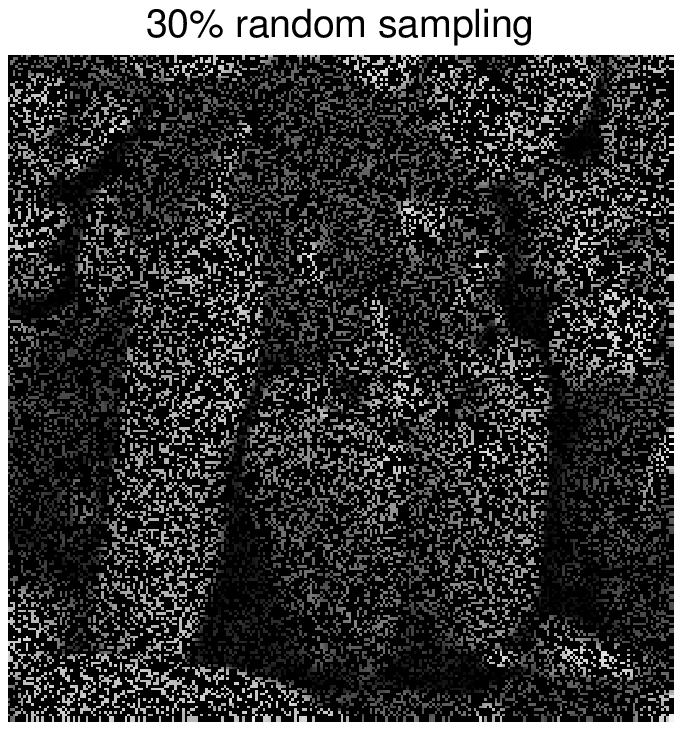}
  \end{minipage}
  \begin{minipage}[t]{0.4\linewidth}
  \centering
  \includegraphics[width=1\textwidth]{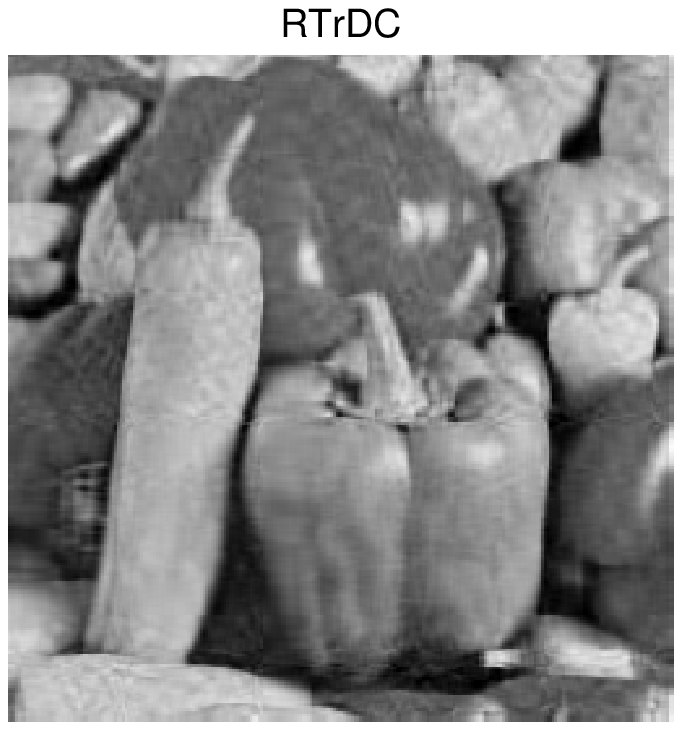}
  \end{minipage}
    \begin{minipage}[t]{0.4\linewidth}
  \centering
  \includegraphics[width=1\textwidth]{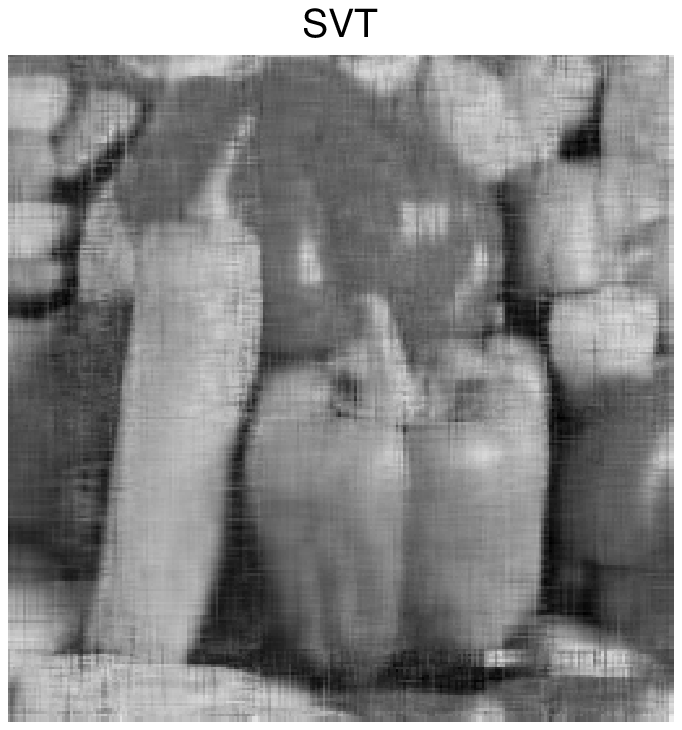}
  \end{minipage}
    \begin{minipage}[t]{0.4\linewidth}
  \centering
  \includegraphics[width=1\textwidth]{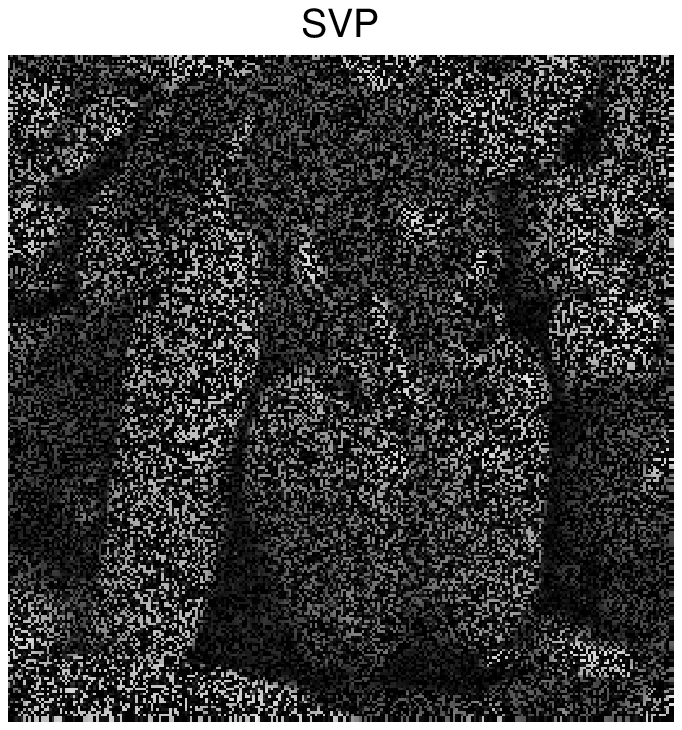}
  \end{minipage}
  \caption{Comparisons of RTrDC algorithm, SVT algorithm and SVP algorithm for recovering the approximated Peppers image (noiseless case) with $\mathrm{SR}=0.30$.} \label{figure4}
\end{figure}

Tables \ref{table1} and \ref{table2} report the numerical results of RTrDC algorithm, SVT algorithm and SVP algorithm for the image inpainting problems with
fixed rank $r=30$. Combined with Figure \ref{figure3} and Figure \ref{figure4}, we can find that our algorithm performs far more better than other two algorithms.

\begin{table}[h!]
\footnotesize
\centering
\begin{tabular}{|c||l|l|l|l|}\hline
\multicolumn{4}{|c|}{Venous image, noiseless, $(r=30,\mathrm{SR}=0.30,\mathrm{FR}= 2.1242)$}\\\hline
Algorithm & RTrDC algorithm& SVT algorithm& SVP algorithm\\
\hline
RE& 9.84e--04 & 2.36e--01 & 8.15e--01 \\
\hline
\multicolumn{4}{|c|}{Peppers image, noiseless, $(r=30,\mathrm{SR}=0.30,\mathrm{FR}=1.3597)$}\\\hline
Algorithm & RTrDC algorithm& SVT algorithm& SVP algorithm\\
\hline
RE & 9.70e--04& 1.04e--01 & 8.13e--01\\
\hline
\end{tabular}
\caption{\scriptsize Numerical results of RTrDC algorithm, SVT algorithm, SVP algorithm for image inpainting problems (noiseless case), $\mathrm{SR}=0.30$.}\label{table2}
\end{table}

\subsection{Image inpainting-noise case }
\begin{figure}[h!]
  \centering
  \begin{minipage}[t]{0.4\linewidth}
  \centering
  \includegraphics[width=1\textwidth]{approximated-venous-30.eps}
  \end{minipage}
  \begin{minipage}[t]{0.4\linewidth}
  \centering
  \includegraphics[width=1\textwidth]{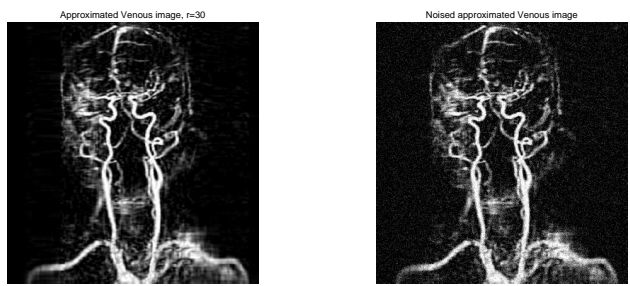}
  \end{minipage}
  \caption{Approximated Venous image and its noised image.} \label{figure5}
\end{figure}

\begin{figure}[h!]
  \centering
  \begin{minipage}[t]{0.4\linewidth}
  \centering
  \includegraphics[width=1.1\textwidth]{approximate-pepper-30.eps}
  \end{minipage}
  \begin{minipage}[t]{0.4\linewidth}
  \centering
  \includegraphics[width=1.1\textwidth]{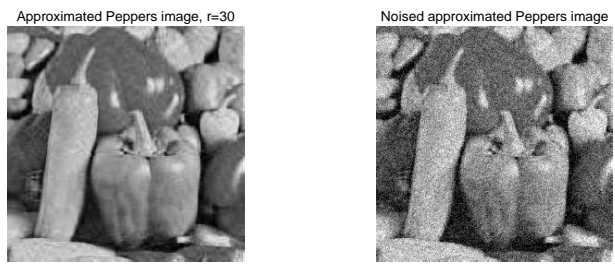}
  \end{minipage}
  \caption{Approximated Peppers image and its noised image.} \label{figure6}
\end{figure}

In this section, we consider the noise case and take a series of experiments to demonstrate the performance of the RTrDC algorithm on two image inpainting problems.
We generate the noised image by
$$\texttt{imnoise(image, 'gaussian', 0, 0.01)}.$$
The approximated Venous image and its noised image are displayed in Figure \ref{figure5}, and the approximated Peppers image and its noised image are displayed in
Figure \ref{figure6}.

\begin{table}[h!]
\footnotesize
\centering
\begin{tabular}{|c||l|l|l|l|}\hline
\multicolumn{4}{|c|}{Venous image, noise $(r=30,\mathrm{SR}=0.40,\mathrm{FR}=2.8323)$}\\\hline
Algorithm & RTrDC algorithm& SVT algorithm& SVP algorithm\\
\hline
RE& 2.02e--01& 3.38e--01 & 8.03e--01\\
\hline
\multicolumn{4}{|c|}{Peppers image, noise, $(r=30,\mathrm{SR}=0.40,\mathrm{FR}=1.8129)$}\\\hline
Algorithm & RTrDC algorithm& SVT algorithm& SVP algorithm\\
\hline
RE & 1.15e--01& 1.62e--01& 7.74e--01 \\
\hline
\end{tabular}
\caption{\scriptsize Numerical results of RTrDC algorithm, SVT algorithm, SVP algorithm for image inpainting problems (noise case), $\mathrm{SR}=0.40$.}\label{table3}
\end{table}

\begin{table}[h!]
\footnotesize
\centering
\begin{tabular}{|c||l|l|l|l|}\hline
\multicolumn{4}{|c|}{Venous image, noise, $(r=30,\mathrm{SR}=0.30,\mathrm{FR}= 2.1242)$}\\\hline
Algorithm & RTrDC algorithm& SVT algorithm& SVP algorithm\\
\hline
RE& 2.04e-01& 3.98e-01 & 8.56e-01\\
\hline
\multicolumn{4}{|c|}{Peppers image, noise, $(r=30,\mathrm{SR}=0.30,\mathrm{FR}=1.3597)$}\\\hline
Algorithm & RTrDC algorithm& SVT algorithm& SVP algorithm\\
\hline
RE & 1.27e-01& 1.92e-01& 8.34e-01\\
\hline
\end{tabular}
\caption{\scriptsize Numerical results of RTrDC algorithm, SVT algorithm, SVP algorithm for image inpainting problems (noise case), $\mathrm{SR}=0.30$.}\label{table4}
\end{table}

\begin{figure}[h!]
  \centering
  \begin{minipage}[t]{0.4\linewidth}
  \centering
  \includegraphics[width=1\textwidth]{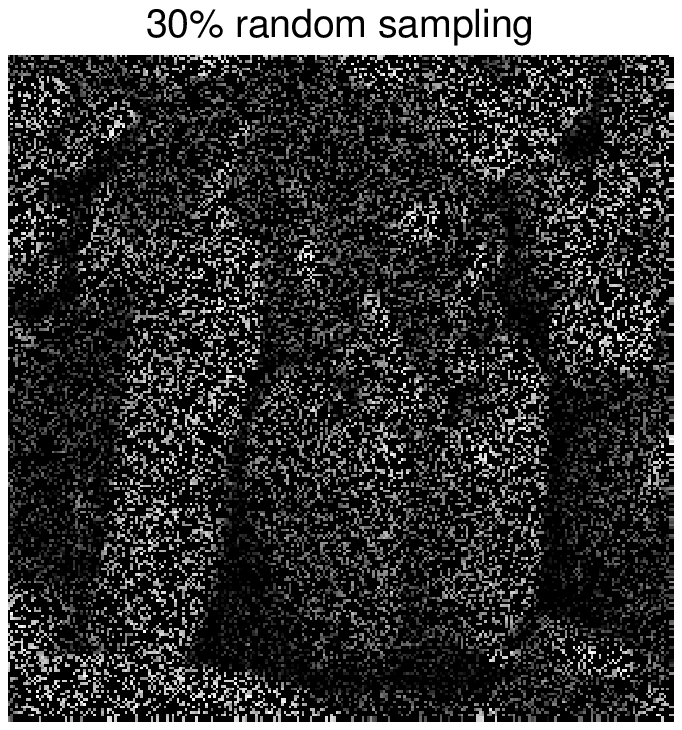}
  \end{minipage}
  \begin{minipage}[t]{0.4\linewidth}
  \centering
  \includegraphics[width=1\textwidth]{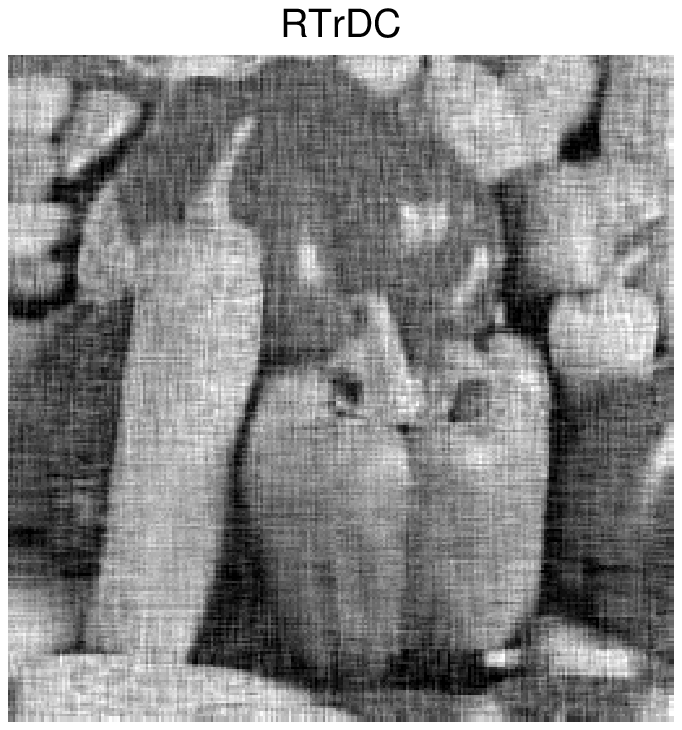}
  \end{minipage}
    \begin{minipage}[t]{0.4\linewidth}
  \centering
  \includegraphics[width=1\textwidth]{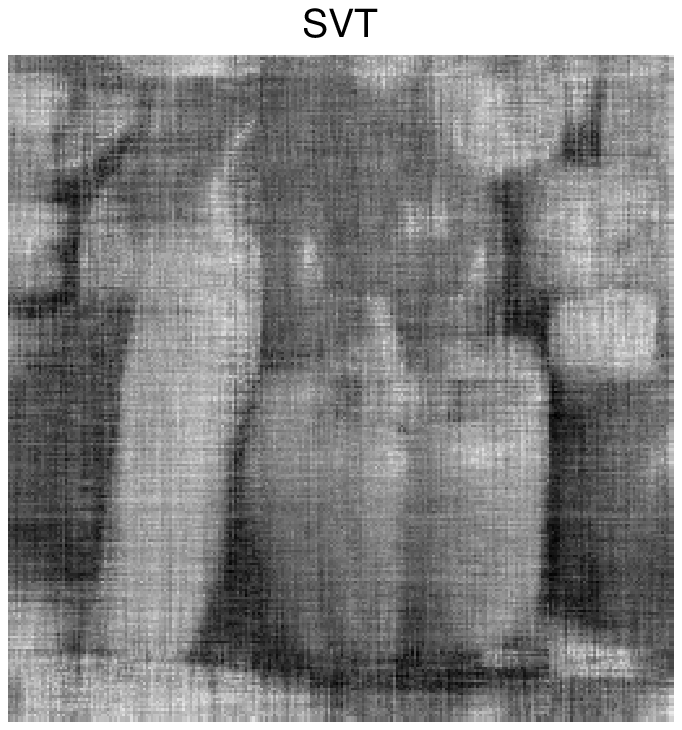}
  \end{minipage}
    \begin{minipage}[t]{0.4\linewidth}
  \centering
  \includegraphics[width=1\textwidth]{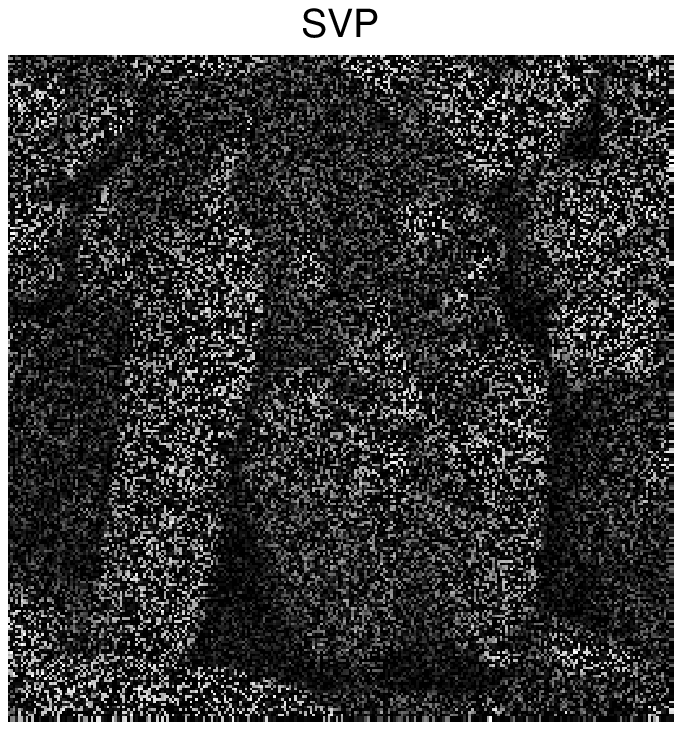}
  \end{minipage}
  \caption{Comparisons of RTrDC algorithm, SVT algorithm and SVP algorithm for recovering the approximated Peppers image (noise case) with $\mathrm{SR}=0.30$.} \label{figure7}
\end{figure}

\begin{figure}[h!]
  \centering
  \begin{minipage}[t]{0.4\linewidth}
  \centering
  \includegraphics[width=1\textwidth]{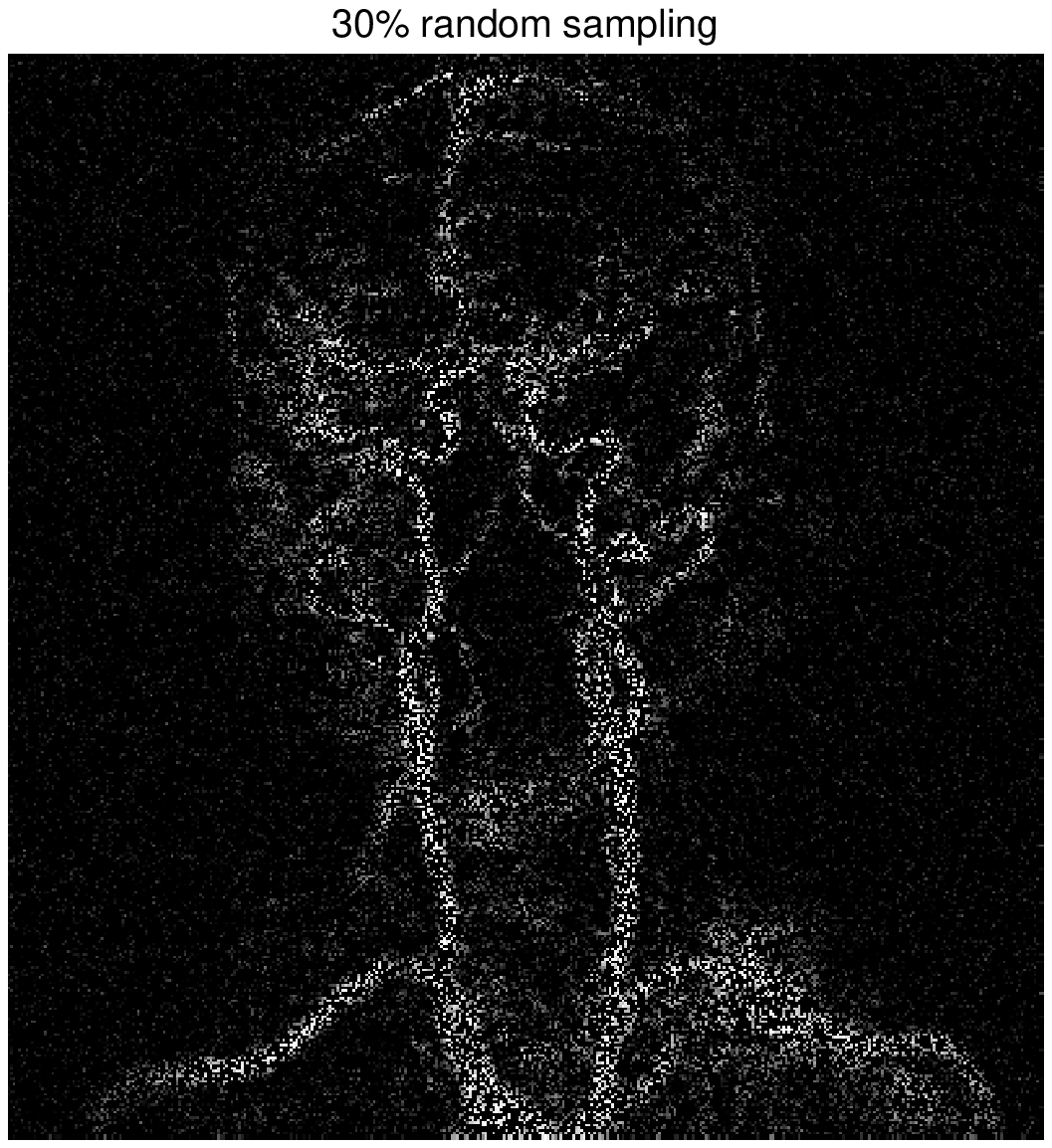}
  \end{minipage}
  \begin{minipage}[t]{0.4\linewidth}
  \centering
  \includegraphics[width=1\textwidth]{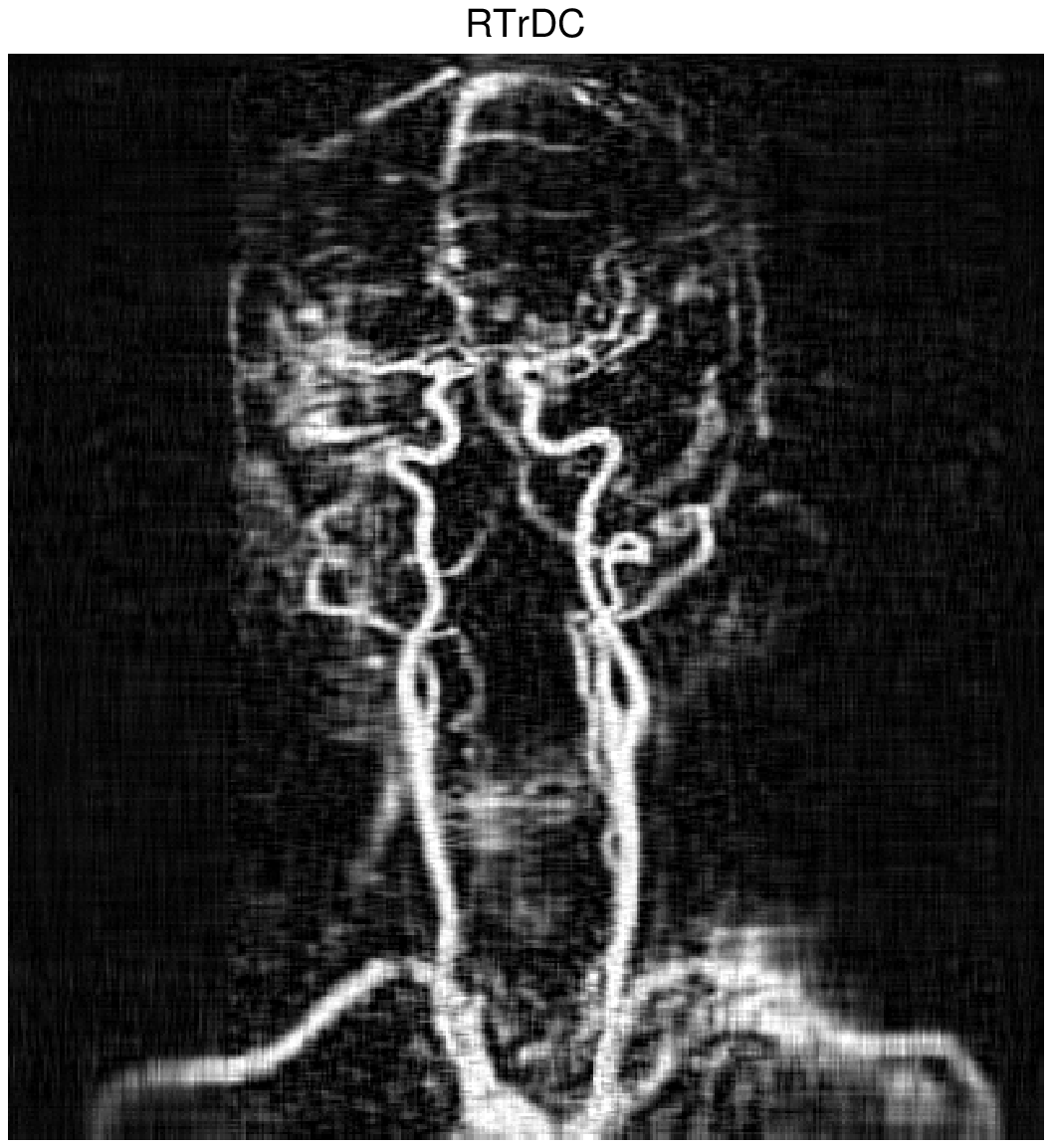}
  \end{minipage}
    \begin{minipage}[t]{0.4\linewidth}
  \centering
  \includegraphics[width=1\textwidth]{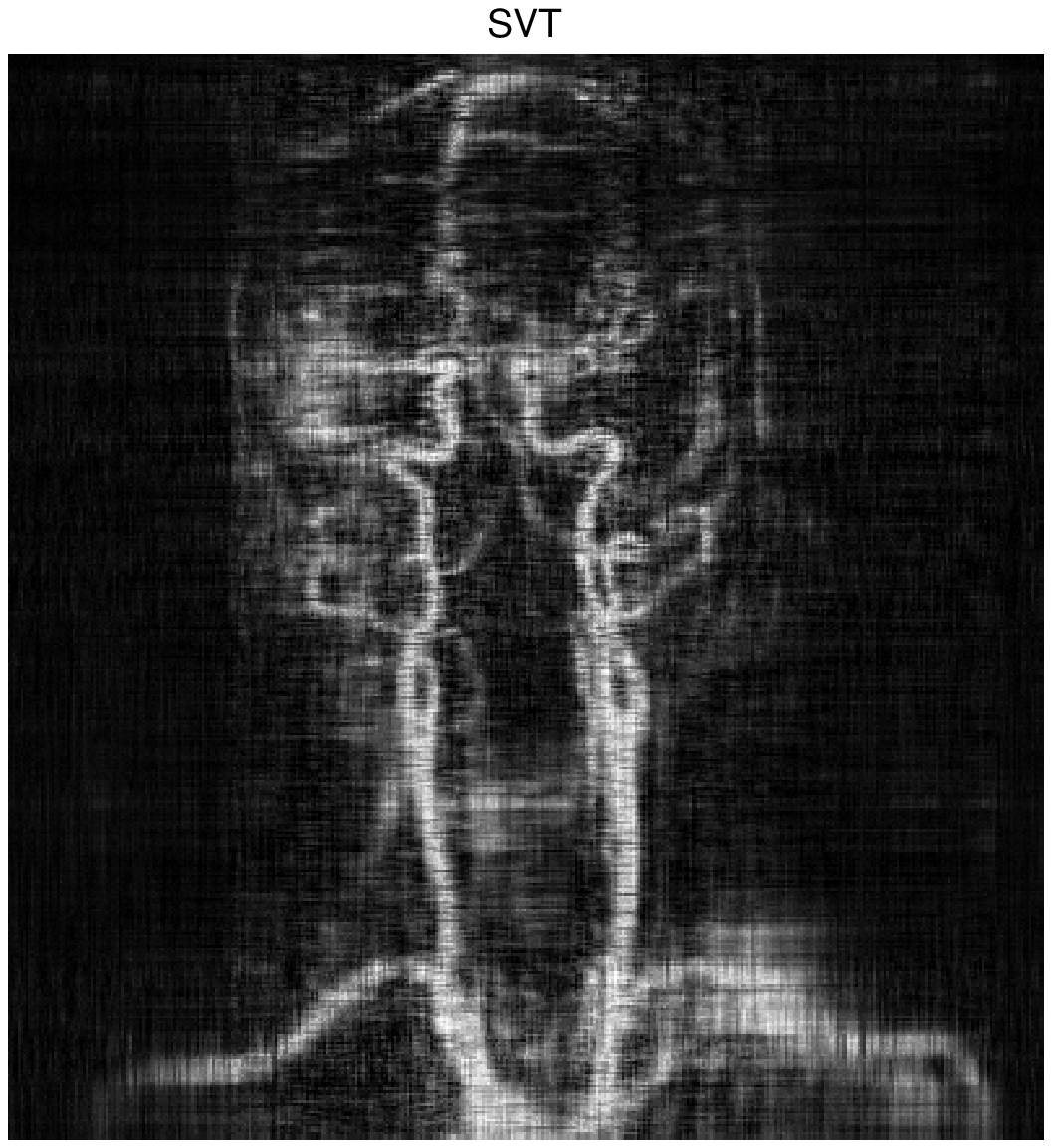}
  \end{minipage}
    \begin{minipage}[t]{0.4\linewidth}
  \centering
  \includegraphics[width=1\textwidth]{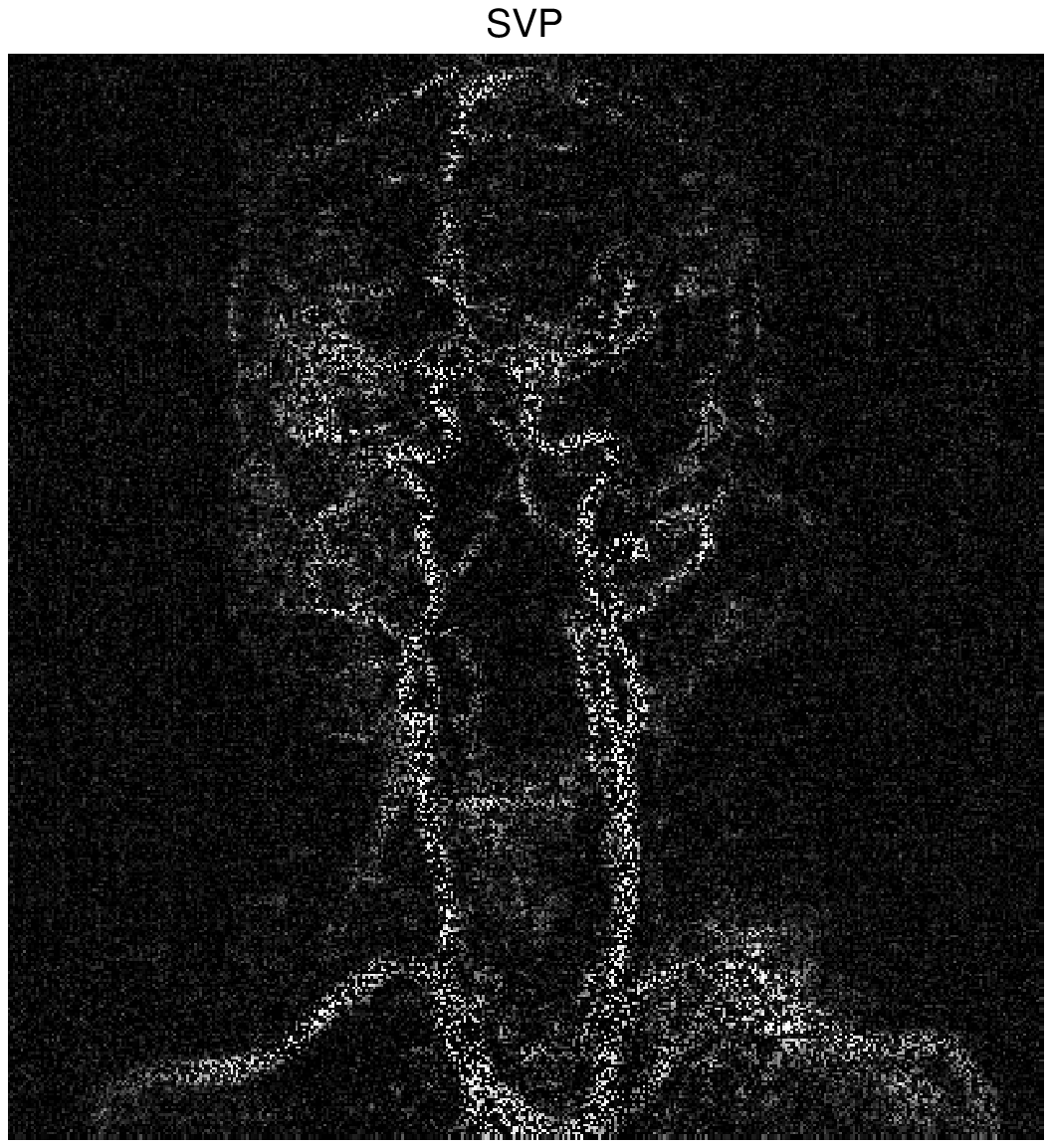}
  \end{minipage}
  \caption{Comparisons of RTrDC algorithm, SVT algorithm and SVP algorithm for recovering the approximated Venous image (noise case) with $\mathrm{SR}=0.30$.} \label{figure8}
\end{figure}

Tables \ref{table3} and \ref{table4} and Figures \ref{figure7} and \ref{figure8} show that the RTrDC algorithm performs the best in finding a low-rank matrix
on image inpainting problems.

\section{Conclusions}\label{section6}
In this paper, a non-convex function is studied to replace the rank function in the problem (AMRM), and translate this NP-hard problem into the problem (TrAMRM).
We theoretically proved that the optimal solution to the problem (TrAMRM) also solves the problem (AMRM) whenever some specific conditions satisfied. Moreover,
we also proved that the optimal solution to the problem (TrAMRM) could be approximately obtained by solving its regularization problem (RTrAMRM) for some proper
smaller $\lambda>0$. Lastly, the DC algorithm is utilized to solve the problem (RTrAMRM). Numerical experiments on image inpainting problems show that our method
performs effectively in recovering low-rank images compared with some art-of-state methods.
\begin{acknowledgements}
We would like to thank editorial and referees for their comments, which may help us to enrich the content and improve the presentation of the results in this paper. The work was supported by the National Natural Science Foundations of China (11771347, 11131006, 41390450, 11761003, 11271297) and the Science Foundations of Shaanxi Province of China (2016JQ1029, 2015JM1012).
\end{acknowledgements}


\begin{thebibliography}{}
%
%


\bibitem{Ji1}
S. Ji, K. F. Sze, Z. Zhou, M. C. So and Y. Ye, Beyond Convex Relaxation: A polynomial-time nonconvex optimization approach to network localization.
IEEE INFOCOM, 2013, 12, 2499-2507 (2013)


\bibitem{Faz2}
M. Fazel, H. Hindi and S. Boyd, A rank minimization heuristic with application to minimum order system approximation. In proceedings of
American Control Conference, Arlington, VA, 6, 4734-4739 (2001)


\bibitem{Faze3}
M. Fazel, H. Hindi and S. Boyd, Log-det heuristic for matrix minimization with applications to Hankel and Euclidean distance matrices.
In Proceedings of American Control Conference, Denever, Colorado, 3, 2156-2162 (2003)


\bibitem{Cand4}
E. J. Cand\`{e}s, B. Recht, Exact matrix completion via convex optimization. Foundations of Computational Mathematics, 9, 717-772 (2009)


\bibitem{Jan5}
D. Jannach, M. Zanker, A. Felfernig and G. Friedrich, Recommender Systerm: An Introduction. Cambridge university press, New York (2012)

\bibitem{Netfix6}
Netfix prize website, https://www.netflixprize.com/


\bibitem{Yeg7}
S. F. Yeganli, R. Yu, Image inpainting via singular value thresholding. The 21st Signal Processing and Communications Applications Conference (SIU), IEEE Conference, Turkey, 1--4 (2013)

\bibitem{Xi8}
X. Peng, C. Lu, Z. Yi and H. Tang, Connections Between Nuclear-Norm and Frobenius-Norm-Based Representations. IEEE Transactions on Neural Networks and Learning Systems, 29(1): 218--224 (2018)

\bibitem{Xi9}
X. Peng, J. Lu, Z. Yi and R. Yan, Automatic Subspace Learning via Principal Coefficients Embedding. IEEE Transactions on Cybernetics, 47(11): 3583--3591 (2017)

\bibitem{Recht10}
B. Recht, M. Fazel and P. A. Parrilo, Guaranteed minimum-rank solution of linear matrix equations via nuclear norm minimization.
SIAM Review, 52(3): 471-501 (2010)

\bibitem{Cand11}
E. J. Candes and T. Tao, The power of convex relaxation: Near-optimal matrix completion. IEEE Transactions on Information Theory, 56 (5): 2053--2080 (2010)

\bibitem{Faze12} 
M. Fazel, Matrix rank minimization with applications. PhD thesis, Stanford University, 2002

\bibitem{Cand13}
E. J. Candes and Y. Plan, Matrix completion with noise. Proceedings of the IEEE, 98(6): 925--936 (2010)

\bibitem{Cai14}
J. Cai, E. J. Candes and Z. Shen, A singular value thresholding algorithm for matrix completion, SIAM Journal on Optimization, 20(4): 1956--1982 (2010)


\bibitem{Ma15}
S. Ma, D. Goldfarb and L. Chen, Fixed point and Bregman iterative methods for matrix rank minimization. SIAM Journal on Optimization, 128(1): 321--353 (2011)


\bibitem{Toh16}
K. C. Toh and S. Yun, An accelerated proximal gradient algorithm for nuclear norm regularized linear least squares problems. Pacific Journal of Optimization, 6(3): 615--640 (2012)



\bibitem{Yu17}
Y. Yu, J. Peng and S. Yue, A new nonconvex approach to low-rank matrix completion with application to image inpainting, https://doi.org/10.1007/s11045-018-0549-5, 2018


\bibitem{Thi18}
L. T. Hoai An and P. D. Tao, The DC (difference of convex functions) Programming and DCA revisited with DC models of real world nonconvex optimization problems. 
Annals of Operations Research, 133(1): 23--46 (205)


\bibitem{Thi19}
R. Horst and N. V. Thoai, DC Programming: Overview. Journal of Optimization Theory and Applications, 103(1): 1--43 (1999)


\bibitem{Lew20}
A.S. Lewis, The convex analysis of unitarily invariant matrix functions. Journal of Convex Analysis, 2(1-2): 173--183 (1995)



\bibitem{Dau21}
I. Daubechies, M. Defrise and C. De Mol, An iterative thresholding algorithm for linear inverse problems with a sparsity constraint. 
Communications on Pure and Applied Mathematics, 57(11): 1413--1457 (2004)


\bibitem{Yu22}
Y. Yu, J. Peng, The Moreau envelope based efficient first-order methods for sparse recovery. Journal of Computational and Applied Mathematics, 322: 109--128 (2017)


\bibitem{Tib23}
R. Tibshirani, Regression Shrinkage and Selection via the Lasso. Journal of the Royal Statistical Society, 58(1): 267--288 (1996)


\bibitem{Jos24}
J. Brodie, I. Daubechies, C. De Mol, D. Giannone and I. Loris, Sparse and stable markowitz portfolio. Proceedings of the National Academy of Sciences of the United States of America, 106(30): 12267--12272 (2009)


\bibitem{Lan25}
X. Lan, A. J. Ma. P. C. Yuen and R. Chellappa, Joint Sparse Representation and Robust Feature-Level Fusion for Multi-Cue Visual Tracking. IEEE Transactions on Image Processing, 24(12): 5826--5841 (2015)


\bibitem{Lan26}
X. Lan, S. Zhang, P. C. Yuen and R. Chellappa, Learning Common and Feature-Specific Patterns: A Novel Multiple-SparseRepresentation-Based Tracker. IEEE Transactions on Image Processing, 27(4): 2022--2037 (2018)


\bibitem{Lan27}
X. Lan, S. Zhang and P. C. Yuen, Robust joint discriminative feature learning for visual tracking. Proceedings of the 25th International Joint Conference on Artificial Intelligence, IJCAI 2016, 3403--3410 (2016)


\bibitem{Lan29}
X. Lan, P. C. Yuen and R. Chellappa, Robust MIL-Based Feature Template Learning for Object Tracking. Proceedings of the Thirty-First AAAI Conference on Artificial Intelligence, AAAI 2017, 4118--4125 (2017)

\bibitem{Lan30}
X. Lan, A. J. Ma and P. C. Yuen, Multi-cue Visual Tracking Using Robust Feature-Level Fusion Based on Joint Sparse Representation. 2014 IEEE Conference on Computer Vision and Pattern Recognition,CVPR 2014, 1194--1201 (2014)



\bibitem{Rag31}
R. Meka, P. Jain and I. S. Dhillon, Guaranteed rank minimization via singular value projection. Proceeding of the Neural Information Processing Systems Conference, NIPS 2009, 937--945 (2009)






































\end{thebibliography}


\end{document}